\documentclass[12pt,fleqn]{article}
\usepackage{latexsym}
\usepackage[dvips]{graphics}
\usepackage{amssymb, amsmath}
\usepackage{bm}

\def\disp{\displaystyle } 
\def\calP{{\cal{P}}}
\def\calS{{\cal{S}}}
\def\calL{{\cal{L}}}
\def\calA{{\cal{A}}}
\def\setR{\mathbb{R}}
\def\setZ{\mathbb{Z}}
\def\veca{\bm{a}}
\def\vecb{\bm{b}}
\def\vecx{\bm{x}}
\def\vecy{\bm{y}}
\def\vecz{\bm{z}}
\def\vect{\bm{t}}
\def\vecn{\bm{n}}
\def\err{\mathrm{err}}
\def\EOC{\mathrm{EOC}}
\def\Image{\mathop\mathrm{Image}}
\def\<{\langle}
\def\>{\rangle}
\def\ratio{\mbox{\sf r}}
\def\eps{\varepsilon}
\makeatletter
 
 \@addtoreset{equation}{section}
\makeatother
\makeatletter
 
 \@addtoreset{figure}{section}
\makeatother
\makeatletter
 
 \@addtoreset{table}{section}
\makeatother
\newtheorem{lm}{Lemma}[section]
\newtheorem{df}[lm]{Definition}
\newtheorem{Q}[lm]{Question}
\newtheorem{ex}[lm]{Example}
\newtheorem{rem}[lm]{Remark}
\begin{document}
\title{\vspace{-2.5truecm}
\begin{flushleft}
{\small{\sf submitted to}: {\it Japan J. Indust. Appl. Math.}}\\
{\small{\sf running title}: Curvature adjusted tangential velocity}\\
\end{flushleft}
\quad\\\large
Evolution of plane curves with a curvature adjusted tangential velocity\thanks{
The authors were/are supported by grants: 
VEGA 1/0657/11 (D\v{S}) and 
Grant-in-Aid for Encouragement of Young Scientists 21740079 (SY).}}
\author{\normalsize
Daniel \v SEV\v COVI\v C\thanks{
Institute of Applied Mathematics, 
Faculty of Mathematics, Physics and Informatics, 
Comenius University, 
842 48 Bratislava, Slovak Republic. 
{\it E-mail}: sevcovic@fmph.uniba.sk}
\ and \ 
Shigetoshi YAZAKI\thanks{
Faculty of Engineering, 
University of Miyazaki, 
1-1 Gakuen Kibanadai Nishi, 
Miyazaki 889-2192, Japan. 
{\it E-mail}: yazaki@cc.miyazaki-u.ac.jp}}
\date{}
\maketitle
\begin{center}
\begin{minipage}[ht]{0.9\linewidth}{\small
{\bf Abstract.}\ 
We study evolution of a closed embedded plane curve with the normal velocity depending on the curvature, 
the orientation and the position of the curve. 
We propose a new method of tangential redistribution of points by curvature adjusted control 
in the tangential motion of evolving curves. 
The tangential velocity distributes grid points along the curve not only uniform 
but also lead to a suitable concentration and/or dispersion depending on the curvature. 
Our study is based on solutions to the governing system of nonlinear parabolic equations 
for the position vector, tangent angle and curvature of a curve. 
We furthermore present a semi-implicit numerical discretization scheme based 
on the flowing finite volume method. 
Several numerical examples illustrating capability of the new tangential redistribution method are 
also presented in this paper. 


}\end{minipage}
\end{center}
\section{\normalsize Introduction}

The purpose of this paper is to study evolution of a plane curve 
$\Gamma^t, t\in[0,T)$, driven by the normal velocity $v$ which 
is assumed to be a function of the curvature $k$, tangential angle $\nu$ 
and position vector $\vecx\in\Gamma^t$, 
\begin{equation}\label{geomrov}
v= \beta(\vecx, \nu, k) \,. 
\end{equation}
Geometric equations of the form (\ref{geomrov}) often arise from various  applied problems like e.g. 
the material science, dynamics of phase boundaries in thermomechanics, 
computational geometry, semiconductors industry, 
image processing and computer vision, etc. 
For a comprehensive overview of industrial applications of the geometric equation having 
the form of (\ref{geomrov}) we refer to a book by Sethian~\cite{Sethian1999}. 

We consider a simple, embedded and closed plane curve $\Gamma$ which is 
parameterized by a smooth function $\vecx:\ \setR/\setZ\supset[0, 1]\to\setR^2$ such that 
$\Gamma=\Image(\vecx)=\{\vecx(u);\ u\in[0, 1]\}$ and 
$|\partial_u\vecx|>0$. 
Hereafter, we denote $\partial_\xi{\sf F}=\partial{\sf F}/\partial\xi$, 
and $|\veca|=\sqrt{\veca\cdot\veca}$ where 
$\veca\cdot\vecb$ is the Euclidean inner product between vectors $\veca$ and $\vecb$. 
The unit tangent vector can be defined as 
$\vect=\partial_u\vecx/|\partial_u\vecx|=\partial_s\vecx$, 
where $s$ is the arc-length parameter and $ds=|\partial_u\vecx|du$. The unit inward normal vector $\vecn$ is defined in such a way that 
$\det(\vect,\vecn)=1$. 
The signed curvature in the direction $\vecn$ is denoted by $k$. 
With this orientation we have $k=\det(\partial_s\vecx, \partial^2_s\vecx)$. 
Let $\nu$ be the angle of $\vect$, i.e., $\vect=(\cos\nu, \sin\nu)^{\mathrm{T}}$ and 
$\vecn=(-\sin\nu, \cos\nu)^{\mathrm{T}}$. 
The problem of evolution of curves is stated as follows: 
For a given initial curve $\Gamma^0=\Image(\vecx^0)$, find a family of planar curves
$\{\Gamma^t\}_{0\leq t<T}$, 
$\Gamma^t=\{\vecx(u, t);\ u\in[0, 1]\subset\setR/\setZ\}$ starting 
from $\Gamma^0=\{ \vecx(u, 0)=\vecx^0(u);\ u\in[0, 1]\}$ and evolving 
in the normal direction according to the velocity given by (\ref{geomrov}). 

As typical examples of $\beta$ one can consider 
$\beta=w(\nu)k$ representing the anisotropic mean curvature flow, $\beta=|k|^{\gamma-1}k$ (power like flow arising in image processing)
and $\beta=w(\vecx, \nu)k+F(\vecx, \nu)$, where 
$w$ is a weight function, $F$ is an external force and $\gamma>0$ is a real parameter. 
The normal velocity $v$ is the normal component of the following evolutionary equation for the position vector $\vecx$: 
\begin{equation}\label{eq:x_t=betaN+alphaT}
\partial_t\vecx=\beta\vecn+\alpha\vect, \quad
\vecx(\cdot, 0)=\vecx^0(\cdot). 
\end{equation}
Here $\alpha$ is the tangential component of the velocity vector. 
Notice that the motion in the tangential direction with a tangential velocity $\alpha$ has no effect of the shape 
of evolving closed curves~\cite[Proposition 2.4]{EpsteinG1987}. 
The shape is determined by the value of the normal velocity $\beta$ only. 
Hence the trivial setting $\alpha\equiv 0$ can be chosen. 
In~\cite{Dziuk1994} Dziuk investigated a numerical scheme for $\beta=k$ and $\alpha=0$. 
However, it was documented by many authors (see e.g.~\cite{HouLS1994, MikulaS1999, MikulaS2001, Yazaki2007} 
and references therein) 
that such a choice of $\alpha$ may lead to various numerical instabilities caused by 
either undesirable concentration and/or extreme dispersion of numerical grid points. 

In order to construct a stable numerical computational scheme, 
various choices of a nontrivial tangential velocity have been proposed and analyzed by many authors. 
We present a brief review of development of nontrivial tangential velocities. 
In~\cite{Kimura1994, Kimura1997} Kimura proposed a uniform redistribution scheme for the case when $\beta=k$ 
by using a special choice of $\alpha$ satisfying the uniform redistribution condition (U): $|\partial_u\vecx|=L^t$, 
where $L^t$ is total length of a curve $\Gamma^t$. 
The author also proved convergence of a numerically discretized solution when assuming uniform distribution of initial grid points. 
In~\cite{HouLS1994}, Hou, Lowengrub and Shelley (see also Hou, Klapper and Si~\cite{HouKS1998}) utilized the condition (U) directly  for the case $\beta=k$. 
More precisely, they derived the form of a tangential velocity 
with $\varphi\equiv 1$ and $\omega\equiv 0$ (see our notation (\ref{eq:tangential_velocity}) 
described in Section~\ref{sec:redistirbution}). 
Such a tangential velocity was also proposed and analyzed independently by Mikula and the first author in~\cite{MikulaS2001}. 
In~\cite[Appendix 2]{HouLS1994}, Hou et al. also pointed out the role of tangential velocity (\ref{eq:tangential_velocity}) 
for a general class of the curvature shape function 
$\varphi(k)$ with $\omega\equiv 0$ as a generalization of condition (U). 
However, there was no systematic explanation of the consequence of such a choice of the tangential velocity 
on grid points redistribution. 
In the present paper, the important role of the so-called curvature adjusted tangential velocity (\ref{eq:tangential_velocity}) 
is emphasized from qualitative and quantitative point of view. 

As we already mentioned, in the paper~\cite{MikulaS2001}, the authors 
derived the expression for the tangential velocity (\ref{eq:tangential_velocity}) with $\varphi\equiv 1$ 
and $\omega\equiv 0$ for a rather general class of planar mean curvature flows satisfying 
the geometric equation (\ref{geomrov}) with $\beta=\beta(\nu, k)$. 
This result can be considered as an improvement of that of~\cite{MikulaS1999} in which 
satisfactory results were obtained only for the case when $\beta=\beta(k)$ is linear or sublinear function 
with respect to $k$. 
Next, in the series of papers~\cite{MikulaS2004a, MikulaS2004b, MikulaS2006}, 
Mikula and the first author proposed a method of asymptotically uniform redistribution. 
In terms of our notation, they derived (\ref{eq:tangential_velocity}) with $\varphi\equiv 1$ 
and nontrivial relaxation function $\omega(t)$ for a general class of normal velocities of the form $\beta=\beta(\vecx, \nu, k)$. 
Their method was also applied to geodesic curvature flows and image segmentation problems. 

Besides these uniform or asymptotically uniform redistribution methods, 
in the so-called crystalline curvature flow, the grid points are distributed densely (sparsely) 
in those part of the curve where the absolute value of the curvature is large (small). 
Although this redistribution is far from being uniform, numerical computation is quite stable. 
One of the reasons for such a behavior is that polygonal curves are restricted to a class of 
admissible facet directions. 
In order to extract essence of the crystalline curvature flow of polygonal curves and generalize it to 
a wide class of plane curve evolution, 
the second author showed that the tangential velocity $\alpha=-\partial_s\beta/k$ 
is implicitly involved in the crystalline curvature flow of planar curves (cf.~\cite{Yazaki2007}). 
Notice that such a tangential velocity is a local one since its value at some point $\vecx$ of a curve 
depends on the local properties of the curve near $\vecx$. 
In this point, it is worthwhile to mention the paper by Deckelnick~\cite{Deckelnick1997} 
who used locally defined tangential velocity having the form $\alpha=-\partial_u(|\partial_u\vecx|^{-1})$ 
for the special case when $\beta=k$. Recently, Pau\v{s} et al. successfully applied curvature adjusted tangential velocity for evolution of 
open curves having important applications in the dislocation dynamics (see~\cite{PausB2009, PausBK2010}). 

The asymptotically uniform redistribution is quite effective and applicable for a wide range of applications. However, from the approximation point of view, unless the solution curve is a circle, there is no reason to accept uniform redistribution for a general case of a curve evolution. 
Hence the desired redistribution should take into account the shape of evolution curves and it should also depend on the modulus of the curvature. 
We will furthermore present a combination of the method of asymptotic uniform redistribution~\cite{MikulaS2004a} and the crystalline tangential velocity discussed in~\cite{Yazaki2007}. 

As far as 3D implementation of tangential redistribution is concerned, 
in the recent paper by Barrett, Garcke and N\"urnberg~\cite{BarrettGN2007}, 
they proposed and studied a new efficient numerical scheme for evolution of surfaces driven by the Laplacian of the mean curvature. 
It turns out that a uniform redistribution of tangential velocity vectors is 
implicitly built in their numerical scheme. 
A similar scheme with an explicit expression for the tangential redistribution has been proposed 
and utilized by Morigi in~\cite{Morigi2010}. 

The organization of the paper is as follows: 
In Section~\ref{sec:GE}, we recall the system of PDEs for geometric equations (\ref{geomrov}) and (\ref{eq:x_t=betaN+alphaT}). 
In Section~\ref{sec:redistirbution}, 
we will address the question how to control redistribution of grid points by using a nontrivial tangential velocity. 
We focus on a class of tangential velocities dependent locally on the curvature combined 
with known asymptotically uniform redistributions. 
The aim of Section~\ref{sec:STATIC} is to further motivate our study by 
constructing a curvature adjusted redistribution yielding the best possible approximation of the evolving family of 
planar curves as far as the minimization of the error between the length of a curve and the length and enclosed of its polygonal approximation is concerned. 
In Section~\ref{sec:scheme}, 
a numerical solution to the system of governing equations will be constructed by means of the so-called flowing finite volume method. 
In Section~\ref{sec:results}, 
we will present various numerical examples of applications of the curvature adjusted tangential velocity. 
Finally, 
we conclude with some key points in Section~\ref{sec:conlcusion}. 
\section{\normalsize Governing equations}\label{sec:GE}

Without loss of generality, we will rewrite the normal velocity 
$\beta(\vecx, \nu, k)$ as follows:
\begin{equation}\label{eq:beta}
\beta=w(\vecx, \nu, k)k+F(\vecx, \nu).
\end{equation}
The tangential velocity $\alpha$ will be defined later in Section~\ref{sec:redistirbution}. 
Using Fren\'et formula $\partial^2_s \vecx =\partial_s\vect = k \vecn$, the following equation for the position vector follows from (\ref{geomrov}) and (\ref{eq:x_t=betaN+alphaT}):
\begin{equation}\label{eq:equation_position}
\partial_t\vecx=w\partial_s^2\vecx+\alpha\partial_s\vecx+F\vecn, 
\end{equation}
where the operators $\partial_s$ and $\partial_s^2$ stand for arc-length derivatives, i.e. $\partial_s{\sf F}=g^{-1}\partial_u{\sf F}$ and 
$\partial_s^2{\sf F}=g^{-1}\partial_u\partial_s{\sf F}$, 
respectively. Here $g=|\partial_u\vecx|$ denotes the local length of a curve parameterized by $\vecx$. 
\begin{rem}{\rm
If $\beta$ is a linear or superlinear function with respect to $k$ like e.g. 
$w(\nu)k$, $w(\vecx, \nu)k+F(\vecx, \nu)$, $|k|^{\gamma-1}k$ $(\gamma\geq 1)$, 
then the function $w$ appearing in (\ref{eq:beta}) has no singularity at $k=0$. 
However, in the sublinear case, a singularity occurs at $k=0$. 
Then a regularization at $k=0$ will be necessary from both theoretical (local existence and uniqueness of solutions) 
as well as practical (construction of a stable numerical approximation scheme) point of view. 
We refer the reader to~\cite{MikulaS2001} for further discussion of this issue. 
For example, if $w=|k|^{\gamma-1}, \ \gamma\in(0, 1)$, 
one can regularize $|k|^{-1}$ as follows: $|\widetilde{k}|^{-1}=|k|^{-1}$ for $|k|>\eps$ and 
$|\widetilde{k}|^{-1}=\eps^{-1}$ for $|k|\leq\eps$, for a small $\eps>0$. 
}\end{rem}

Following~\cite{MikulaS2001}, one can derive a closed system of PDEs governing the motion of 
curves satisfying the geometric equation (\ref{geomrov}) and (\ref{eq:x_t=betaN+alphaT}). These equations can be derived using the facts that $k=\det(\partial_s\vecx, \partial_s^2\vecx)$, 
$\vect=\partial_s\vecx=(\cos\nu, \sin\nu)^{\mathrm{T}}$, Frenet's formulae $\partial_s\vect=k\vecn$, $\partial_s\vecn=-k\vect$ and 
the relation $k=\partial_s\nu$. The resulting system of governing PDEs together with (\ref{eq:equation_position}) reads as follows:
\begin{eqnarray}\disp
&& \partial_tg=\left(-k\beta+\partial_s\alpha\right)g, 
\label{eq:local_length}\\[5pt]\disp
&& \partial_tk=\partial_s^2\beta+\alpha\partial_sk+k^2\beta, 
\label{eq:equation_curvature}\\[5pt]\disp
&& \partial_t\nu=(\partial_k\beta)\partial_s^2\nu+(\alpha+\partial_\nu\beta)\partial_s\nu+\nabla_{\vecx}\beta.\vect, 
\label{eq:tangent_angle}
\end{eqnarray}
for $u\in[0, 1]$ and $t\in(0, T)$, with initial conditions 
$g(\cdot, 0)=g^0(\cdot)$, 
$k(\cdot, 0)=k^0(\cdot)$, 
$\nu(\cdot, 0)=\nu^0(\cdot)$, 
$\vecx(\cdot, 0)=\vecx^0(\cdot)$, 
and periodic boundary conditions for $g, k, \vecx$ for $u\in[0,1]\subset\setR/\setZ$ and 
the boundary condition $\nu(1, t)=\nu(0, t)+2\pi$. 
The partial derivative $\nabla_{\vecx}\beta$ is defined as 
$\nabla_{\vecx}\beta=(\partial_{x_1}\beta(\vecx, \cdot, \cdot), \partial_{x_2}\beta(\vecx, \cdot, \cdot))^{\mathrm{T}}$ 
for $\vecx=(x_1, x_2)^{\mathrm{T}}$. 
The initial functions $g^0, k^0, \vect^0$ and $\vecx^0$  are required to satisfy the following compatibility conditions: 
$g^0=|\partial_u\vecx^0|>0$, $\partial_u\vect^0=g^0k^0$, 
$k^0=(g^0)^{-3}\det(\partial_u\vecx^0, \partial_u^2\vecx^0)$. 

In our simple and fast numerical approximation scheme proposed in Section~\ref{sec:scheme} 
we will discretize (\ref{eq:equation_position}) and the tangential velocity equation 
(\ref{eq:tangential_velocity}) with the constraint (\ref{eq:tangential_velocity_av=0}) only. 
With regard to \cite{MikulaS2001}, in the case where a more accurate numerical scheme and results are required, it is recommended to discretize the full system of evolutionary 
PDEs (\ref{eq:equation_position}), (\ref{eq:local_length}), (\ref{eq:equation_curvature}) and 
(\ref{eq:tangent_angle}) for all remaining geometric quantities $g$, $k$, and $\nu$. 
On the other hand, 
as far as the time complexity of computation and simplicity of a numerical approximation scheme are concerned, 
discretization of (\ref{eq:equation_position}) together with (\ref{eq:tangential_velocity}), 
(\ref{eq:tangential_velocity_av=0}), yields satisfactory numerical results.

However, as $w$ may depend on $k$ and $\nu$, 
the proof of existence and uniqueness of a smooth solution to the system of equations 
(\ref{eq:equation_position}), (\ref{eq:tangential_velocity}), (\ref{eq:tangential_velocity_av=0}) 
does not seem to be a straightforward issue. 

Recall that in the case where $\alpha$ is given by the expression (\ref{eq:tangential_velocity}) with $\varphi\equiv 1$ and $\omega=0$, 
the short time existence and uniqueness of smooth solutions to the system of PDEs 
(\ref{eq:equation_position}), (\ref{eq:local_length}), (\ref{eq:equation_curvature}) and (\ref{eq:tangent_angle}) 
has been shown in~\cite{MikulaS2001}. In the forthcoming paper~\cite{SevcovicY20xx} the authors proved local existence and uniqueness of 
a classical solution to the full system of governing equations 
(\ref{eq:equation_position})--(\ref{eq:tangent_angle}) and (\ref{eq:tangential_velocity}), (\ref{eq:tangential_velocity_av=0}).

\section{\normalsize A curvature adjusted tangential redistribution of grid points}\label{sec:redistirbution}

A key tool for construction of a suitable tangential redistribution functional $\alpha$ is the so-called relative local length quantity (cf.~\cite{MikulaS2001}). It is defined by the following ratio $\ratio$: 
\[
\ratio(u, t)=\frac{g(u, t)}{L^t}, \quad
u\in[0, 1], \quad t\in[0, T), 
\]
where $L^t$ is the total length of $\Gamma^t$:
\[
L^t=\int_{\Gamma^t}ds=\int_0^1g\,du, \quad t\in[0, T), 
\]
and $T>0$ is the maximal time of existence of a solution. 

The role of $\ratio$ can be explained as follows: 
suppose that $N$ grid points $\{\vecx(u_i, t)\}_{i=1}^N$ are distributed on the curve $\Gamma^t$ 
for $u_i=i/N$, $i=1, 2, \cdots, N$. Since the arc-length $s$ is given by $\disp s(u, t)=s(0, t)+\int_0^ug(u, t)\,du$, we have 
\[
s(u_{i}, t)-s(u_{i-1}, t)
=\int_{u_{i-1}}^{u_i}g(u, t)\,du
=L^t\int_{u_{i-1}}^{u_i}\ratio(u, t)\,du
\]
for each $i$. Hence, if $\ratio(u, t)\equiv 1$ holds for all $u$ at a time $t$, then all grid points are distributed uniformly in the sense that 
$s(u_{i}, t)-s(u_{i-1}, t)=L^t/N$ for each $i$ at the time $t$. Furthermore, if the evolving family of curves fulfills the limit  $\ratio(u, t)\to 1$ for all $u$ as $t\to T$ then redistribution of grid points become asymptotically uniform as $t\to T$ 
(cf.~\cite{MikulaS2004a, MikulaS2004b}). 

It has been documented by many practical examples (see e.g.~\cite{MikulaS2004a, MikulaS2004b}) 
that (asymptotically) uniform redistribution can significantly stabilize and 
speed up numerical computation. 
We recall that the tangential motion has no influence on the shape of evolved closed planar curves. 
A natural question arises:
\begin{Q}{\rm \label{Q:Question}
How to design the relative local length ratio $\ratio$ and, subsequently, 
$\alpha$ such that redistribution takes into account the shape of the limiting curve. 
In other words, how to densely (sparsely) redistribute grid points on those parts of a curve where 
the modulus of the curvature is large (small). 
}\end{Q}

The modulus $|k|$ of curvature will be measured by the following shape function:
\begin{df}{\rm
Let $\varphi: \setR\to\setR_+$ be a nonnegative even function 
$\varphi(-k)=\varphi(k)>0$ for $k\not=0$ such that $\varphi(k)$ is nondecreasing for $k>0$. 
}\end{df}

As an example of a shape function one can consider
$\varphi(k)\equiv 1$, 
$\varphi(k)=|k|$, 
$\varphi(k)=1-\eps+\eps|k|$ ($\eps\in[0, 1]$), 
$\varphi(k)=\sqrt{\eps^2+k^2}$ ($|\eps|\ll 1$), etc. 
In numerical computations contained in this paper, a linear combination of these functions will be often used (see Example~\ref{ex:varphi=combination}). 

In order to answer Question~\ref{Q:Question}, we introduce a generalized relative local length adopted to the shape function $\varphi$ as follows:
\begin{equation}
\ratio_{\varphi}(u, t)=\frac{g(u, t)}{L^t}\frac{\varphi(k(u, t))}{\<\varphi(k(\cdot, t))\>}, \quad
u\in[0, 1], \quad t\in[0, T). 
\label{rellength}
\end{equation}
Here the bracket function $\<{\sf F}\>$ denotes the average of function ${\sf F}(u, t)$ along the curve $\Gamma^t$: 
\[
\<{\sf F}(\cdot, t)\>
=\frac{1}{L^t}\int_{\Gamma^t}{\sf F}\,ds
=\frac{1}{L^t}\int_{0}^1{\sf F}(u, t)g(u, t)\,du. 
\]
Next we explain the role of the generalized ratio $\ratio_{\varphi}$. 
Suppose, for a moment, that $\ratio_{\varphi}(u, t)\equiv 1$ for all $u$ at a time $t$. 
Then $s(u_{i}, t)-s(u_{i-1}, t)\lessgtr L^t/N$ holds, if $|k|$ satisfies $\varphi(k)\gtrless\<\varphi(k)\>$ 
on the $i$-th interval $[u_{i-1}, u_i]$. 
Using this property we are in a position to construct 
the desired curvature adjusted tangential redistribution which takes account of 
deviations of the shape function $\varphi(k)$ from its average value $\<\varphi(k)\>$. 
Indeed, suppose that we are able to construct a tangential velocity $\alpha$ in such a way that 
$\ratio_{\varphi}(u, t)\to 1$, i.e., 
\[
\theta(u, t)=\ln\ratio_{\varphi}(u, t)\to 0
\]
holds for all $u\in[0, 1]$ as $t\to T$. 
For $t$ close to $T$, we can conclude if $|k|$ is above/below the average in the sense that 
$\varphi(k)\gtrless\<\varphi(k)\>$, 
then $g\lessgtr L^t$ holds, respectively. 
Distribution of grid points on corresponding subarcs is dense/sparse, respectively. 
\begin{ex}{\rm
Consider $\varphi(k)\equiv 1$. Then $\theta=\ln\ratio$. If one constructs  $\alpha$ in such a way that $\theta(u, t)\equiv\theta(u, 0)$ for all $u$ and $t$ then relative local length $\ratio$ is preserved. On the other hand, if one takes $\alpha$ such that $\theta(u, t)\to 0$ (i.e., $\ratio\to 1$) for all $u$ as $t\to T$ then redistribution becomes asymptotically uniform. See Example~\ref{ex:varphi=1}. 
}\end{ex}

Convergence 
$\lim_{t\to T}\ratio_{\varphi}(u, t)\to 1$ is fulfilled provided that there exists a relaxation function $\omega\in L_{\rm loc}^1[0, T)$ such that
\[
\int_0^T\omega(\tau)\,d\tau=\infty
\]
and $\theta(u, t)=\ln((e^{\theta(u, 0)}-1)e^{-\int_0^t\omega(\tau)\,d\tau}+1)$ 
for all $u\in[0, 1]$ and $t\in[0, T)$. 
The previous equation can be rewritten as an ODE:
\begin{equation}\label{eq:ratio_equation}
\partial_t\theta(u, t)+\omega(t)(1-e^{-\theta})=0. 
\end{equation}
\begin{rem}{\rm
For $\tilde{\theta}=\ln\ratio_{\varphi}$ 
one can obtain another convergence such as 
$\tilde{\theta}(u, t)=\tilde{\theta}(u, 0)e^{-\int_0^t\omega(\tau)\,d\tau}$, 
then we have the equivalent ODE $\partial_t\tilde{\theta}(u, t)+\omega(t)\tilde{\theta}(u, t)=0$. 
The convergence rate of $\tilde{\theta}(u, t)$ is the same as that of $\theta(u, t)$. 
}\end{rem}

With regard to~\cite{MikulaS2004a, MikulaS2004b} one can choose the relaxation function $\omega(t)$ as follows:
\[
\omega(t)=\kappa_1-\kappa_2\partial_t\ln L^t
=\kappa_1+\kappa_2\<k\beta\>, 
\]
where $\kappa_1\geq 0$ and $\kappa_2\geq 0$ are nonnegative constants. 
Here we have employed the total length equation:
$\partial_t L^t+\int_{\Gamma^t}k\beta\,ds=0$ (cf.~\cite{MikulaS2001}). 
If we formally set $\kappa_1=\kappa_2=0$, then the function $\theta$ is constant in time $t$. 
On the other hand, if the maximal existence time is infinite $T=\infty$, 
we can choose $\kappa_1>0$ and $\kappa_2=0$. 
If $T<\infty$ and $L^t\to 0$ as $t\to T$, we can choose $\kappa_2>0$. 
In both cases we have $\int_0^T\omega(\tau)\,d\tau =\infty$ and so $\ratio_{\varphi}(u,t)\to 1$ as $t\to T$. 

The equation for the tangential velocity $\alpha$ can be derived as follows. 
Differentiating $\theta=\ln\ratio_{\varphi}=\ln g+\ln\varphi-\ln\int_0^1\varphi g\,du$ with respect to $t$, 
using (\ref{eq:local_length}), (\ref{eq:equation_curvature}) and the relation 
$\partial_s(\varphi\alpha)=\alpha(\partial_sk)\varphi'+(\partial_s\alpha)\varphi$, 
we obtain 
\[
\begin{split}
\varphi\partial_t\theta
&=(-k\beta+\partial_s\alpha)\varphi
+(\partial_s^2\beta+k^2\beta+\alpha\partial_sk)\varphi'
-\frac{\varphi}{L^t\<\varphi\>}\int_0^1\partial_t(\varphi g)\,du \\
&=\partial_s(\varphi\alpha)-f+\frac{\varphi}{\<\varphi\>}\<f\>, 
\end{split}
\]
where $f=\varphi(k)k\beta-\varphi'(k)\left(\partial_s^2\beta+k^2\beta\right)$, $
\varphi'(k)=\partial_k\varphi(k)$. 
Hence (\ref{eq:ratio_equation}) holds 
if and only if the tangential velocity $\alpha$ satisfies the following equation:
\begin{equation}\label{eq:tangential_velocity}
\frac{\partial_s(\varphi(k)\alpha)}{\varphi(k)}
=\frac{f}{\varphi(k)}-\frac{\<f\>}{\<\varphi(k)\>}
+\omega(t)\left(\ratio_{\varphi}^{-1}-1\right). 
\end{equation}

Equation (\ref{eq:tangential_velocity}) is written in the form 
$\partial_s(\varphi(k)\alpha)={\mathcal F}$, where 
\[
{\mathcal F}:=
f-\frac{\langle f\rangle}{\langle\varphi(k)\rangle}\varphi(k)
+\omega(t)\varphi(k)(\ratio_\varphi^{-1}-1). 
\]
Since 
$\int_{\Gamma^t}\varphi(k)(\ratio_\varphi^{-1}-1)\,ds
=\int_0^1(L^t\langle\varphi(k)\rangle g^{-1}-\varphi(k))g\,du=0$, 
we obtain $\langle{\mathcal F}\rangle=0$. 
Thus the equation $\partial_s(\varphi(k)\alpha)={\mathcal F}$ has at least one solution $\alpha$. 
In order to construct a unique solution $\alpha$, we assume the following renormalization condition for $\alpha$:
\begin{equation}\label{eq:tangential_velocity_av=0}
\<\varphi(k)\alpha\>=0. 
\end{equation}

\begin{ex}[uniform redistribution]{\rm \label{ex:varphi=1}
In the case where $\varphi(k)\equiv 1$, 
the tangential velocity equation becomes 
\[
\partial_s\alpha=k\beta-\<k\beta\>+\omega(t)(\ratio^{-1}-1). 
\]
If we set $\omega\equiv 0$, then a solution $\alpha$ preserves the relative local length $\ratio$ (see~\cite{MikulaS2001}). 
Under the assumption $\omega\not\equiv 0$ with suitable $\kappa_1,\kappa_2$, 
redistribution of grid points becomes asymptotically uniform~\cite{MikulaS2004a, MikulaS2004b, MikulaS2006}. 
}\end{ex}
\begin{ex}[crystalline tangential velocity]{\rm \label{ex:varphi=k}
Suppose that the evolving curve $\Gamma^t$ is convex. 
If we consider the shape function $\varphi(k)=|k|$ and $\omega(t)\equiv 0$, 
then, with regard to (\ref{eq:tangential_velocity}) 
we have $\partial_s(k\alpha)=-\partial_s^2\beta$. 
Taking into account renormalization constraint $\<\varphi(k)\alpha\>=0$ 
we end up with $\alpha=-\partial_s\beta/k$. 
This is exactly the same  tangential velocity as it was derived by the second author 
in the continuous limit of the crystalline curvature flow (see~\cite{Yazaki2007}). 
If the evolving curve $\Gamma^t$ has negative curvature on those parts, 
the shape function $\varphi(k)=|k|$ is regularized such as, 
for instance, $\varphi(k)=\sqrt{1-\eps+\eps k^2}$ with $\eps\in(0, 1)$. 
}\end{ex}
\begin{ex}[their linear combination]{\rm \label{ex:varphi=combination}
In our numerical computations in Section~\ref{sec:results}, 
we will use the following smoothed shape function  $\varphi$:
\[
\varphi(k)=1-\eps+\eps\sqrt{1-\eps+\eps k^2}, \quad
\eps\in(0, 1). 
\]
It is a linear combination of shape functions in Examples~\ref{ex:varphi=1} and~\ref{ex:varphi=k}. 
Notice that $\varphi(k)\to 1$ if $\eps\to 0$, $\varphi(k)\to |k|$ if $\eps\to 1$, and $\varphi(k)\geq\varphi(0)>0$ for $\eps\in(0, 1)$. 

In Figure~\ref{fig:distribution} we plot an example of distribution of grid points along the ellipse. 
Let $\overline{\vecx}(l)=(a\cos(2\pi l), b\sin(2\pi l))^{\mathrm T}$ be the ellipse with axes ratio 
$a=3:b=1$, where $l\in[0, 1]\subset\setR/\setZ$. 
We will construct a new parameterization $u\in[0, 1]$ using the reparameterization function $l(u)$ such that $\vecx(u)=\overline{\vecx}(l(u))$, $l(0)=0$, $l(1)=1$, $\partial_ul(u)>0$ as follows: since $\partial_u\vecx(u)=\partial_l\overline{\vecx}(l(u))\partial_ul(u)$,  we have $g(u)=\overline{g}(l(u))\partial_ul(u)$, where  $\overline{g}(l)=|\partial_l\overline{\vecx}(l)|$.  Then if $\ratio_{\varphi}\equiv 1$, we obtain the ODE
\[
\partial_ul(u)=\frac{L\<\varphi(\overline{k})\>}{\overline{g}(l(u))\varphi(\overline{k}(l(u)))}
\]
to be solved. We have used the relation $k(u)=\overline{k}(l(u))$. 
Applying the Runge-Kutta ODE solver we calculated results depicted in Figure~\ref{fig:distribution}~(a), (b) and (c). 
Crystalline curvature redistribution $\{\overline{\vecx}(l_i)\}$ is obtained from 
the condition $\overline{\vect}(l_i)=(-\sin(2\pi i/N), \cos(2\pi i/N))^{\mathrm T}$, 
and the $i$-th inward normal vector of the circumscribed polygon is $\overline{\vecn}(l_i)$ 
for $i=1, 2, \cdots, N$ (Figure~\ref{fig:distribution}~(d)). 
}\end{ex}
\section{\normalsize Optimal redistribution of points for stationary curves in the plane}\label{sec:STATIC}

The aim of this section is to further motivate the study of curvature adjusted tangential redistribution. 
In what follows, we will address the question on what is the optimal redistribution of vertices 
$\{\vecx_1, \vecx_2, \cdots, \vecx_N\}$ belonging to a given smooth closed curve $\Gamma$ such that 
the discrepancy between the length/the area of $\Gamma$ and that of a polygon spanned by vertices 
$\{\vecx_1, \vecx_2, \cdots, \vecx_N\}$ is minimal. 
Clearly, in the case $\Gamma$ being a circle, the optimal redistribution is represented by a regular $N$-polygon. 
However, it should be obvious that, e.g. for an oval, 
the optimal redistribution has to take curvature of the curve into account. 
Interestingly enough, we will prove that the length/the area discrepancy minimizing redistribution 
$X=\{\vecx_1, \vecx_2, \cdots, \vecx_N\}$, for $N\to\infty$, 
is closely related to the curvature adjusted redistribution discussed in previous section 
with the shape functions 
\begin{align*}
&\varphi(k)=|k|^{2/3}\quad\mbox{for the length discrepancy}, \\
&\varphi(k)=|k|^{1/3}\quad\mbox{for the area discrepancy}. 
\end{align*}
Suppose that a given smooth closed curve $\Gamma$ is parameterized by a smooth function 
$\vecx:[0,1]\mapsto {\mathbb R}^2$. 
Denote by $X=\{\vecx_1, \vecx_2, \cdots, \vecx_N\}$ the set of grid points of $\Gamma$ such that 
$\vecx_i=\vecx(u_i)$ where $u_i=ih$ and $h=1/N$. 

{\bf The length discrepancy.}\ 
Let $L=\int_\Gamma ds$ be the length of $\Gamma$ and 
\[
{\mathcal L}(X)=\sum_{j=0}^{N-1}\vert\vecx_{j+1}- \vecx_j\vert
\]
be the length of a polygon with vertices 
$\{\vecx_1, \vecx_2, \cdots, \vecx_N\}$. 
Here we have identified $\vecx_{N+i}=\vecx_i$ since $\Gamma$ is assumed to be a closed curve. 
Clearly, ${\mathcal L}(X) \le L$. 
Our goal is to find conditions under which the parameterization $\vecx(\cdot)$ yields the minimizer 
$X=\{\vecx_1, \vecx_2, \cdots, \vecx_N\}$ of the problem
\begin{equation}
\min_{X\subset\Gamma} (L - {\mathcal L}(X)). 
\label{minimization}
\end{equation}
Since 
\[
{\mathcal L}(X) = 
\sum_{j=0}^{i-2} \vert \vecx_{j+1}- \vecx_j\vert
+
\vert \vecx_{i}- \vecx_{i-1}\vert
+
\vert \vecx_{i+1}- \vecx_{i}\vert
+
\sum_{j=i+1}^{N-1} \vert \vecx_{j+1}- \vecx_j\vert
\]
we obtain the following expression for the derivative of ${\mathcal L}$ with respect to 
$\vecx_i$ in the direction $\vecy\in {\mathbb R}^2$:
\[
{\mathcal L}^\prime_{\vecx_i}(X) \vecy = 
-\widetilde{\vecn}_i.\vecy, \quad
\widetilde{\vecn}_i
=\frac{\vecx_{i+1}-\vecx_i}{\vert \vecx_{i+1}-\vecx_i\vert}
-\frac{\vecx_i-\vecx_{i-1}}{\vert \vecx_i-\vecx_{i-1}\vert}
\]
Recall that the above minimization problem (\ref{minimization}) is subject to the constraint 
$\vecx_i \in \Gamma$ for each $i=1, 2, \cdots, N$. 
Therefore, the first order necessary condition for the constrained minimizer $X$ of (\ref{minimization}) 
reads as follows:
\[
{\mathcal L}^\prime_{\vecx_i}(X) \vecy=0\quad 
\hbox{for any} \ \vecy \sim \vect_i,\quad 
i=1, 2, \cdots, N, 
\]
where the symbol $\vecy\sim \vect_i$ means that the vector $\vecy$ is collinear with 
the unit tangent vector $\vect_i$ to $\Gamma$ at the point $\vecx_i\in\Gamma$. 
Hence, the necessary condition for a set $X=\{\vecx_1, \cdots, \vecx_N\}\in {\mathbb R}^{2\times N}$ of points 
belonging to $\Gamma$ to be a minimizer of the functional $X\mapsto L-{\mathcal L}(X)$ can be rewritten as:
\begin{equation}
\widetilde{\vecn}_i\perp\vect_i,\quad 
i=1, 2, \cdots, N. 
\label{necessary-condL}
\end{equation}
A graphical description of the necessary condition is depicted in Figure~\ref{fig:necessary-cond}~(a). 

Next we will express the necessary condition (\ref{necessary-condL}) 
in terms of the parameterization $\vecx(\cdot)$ of the curve $\Gamma$. 
The Taylor expansion of $\vecx(\cdot)$ at $\vecx_i=\vecx(i h)$ yields
\[
\vecx_{i\pm1}- \vecx_i = \pm \partial_u \vecx(u_i) h 
+\partial^2_u \vecx(u_i) \frac{h^2}{2}
\pm\partial^3_u \vecx(u_i) \frac{h^3}{6}
+\partial^4_u \vecx(u_i) \frac{h^4}{24}  + O(h^5), 
\]
as $h\to 0^+$. Derivatives $\partial^j_u \vecx$ at $u_i = i h$ can be decomposed as follows:
\[
\partial^j_u \vecx(u_i) = b_j \vecn_i  + a_j \vect_i, \quad j=1,2, \cdots, 
\]
where $\vecn_i$ and $\vect_i$ are, respectively, 
the unit inward normal and tangent vector to the curve $\Gamma$ at a point $\vecx_i\in\Gamma$. 
Then $a_1= g$, $b_1=0$ where $g=g_i=\vert \partial_u \vecx(u_i)\vert$ 
is the local length at the point $\vecx_i\in\Gamma$. 
Furthermore, using the Fren\'et formulae $g^{-1}\partial_u \vect= \partial_s \vect = k \vecn$ and 
$g^{-1}\partial_u \vecn= \partial_s \vecn = - k \vect$, 
we obtain the recurrent relations
\begin{equation}
a_1 = g, \quad b_1 =0, \quad
a_{j+1} = \partial_u a_j - g k b_j,\quad 
b_{j+1} = \partial_u b_j + g k a_j.
\label{recurrent}
\end{equation}
for $j\ge 1$ where $k=k_i$ is the curvature of $\Gamma$ at $\vecx_i\in\Gamma$.

Let us define $\Phi(h)$ and $\Psi(h)$ as the following auxiliary functions
\[
\Phi(h)= \sum_{j=1}^\infty \frac{a_j}{j!} h^{j-1}, 
\quad 
\Psi(h)= \sum_{j=1}^\infty \frac{b_j}{j!} h^{j-1}. 
\]
Using these functions, the forward and backward difference at $\vecx_i$ can be expressed as
\[
\frac{\vecx_{i\pm 1}-\vecx_i}{\pm h}
=\Phi(\pm h)\vect_i+\Psi(\pm h)\vecn_i, \quad
\frac{|\vecx_{i\pm 1}-\vecx_i|}{h}
=\sqrt{\Phi(\pm h)^2+\Psi(\pm h)^2}. 
\]
Then the necessary condition (\ref{necessary-condL}) can be rewritten as:
\[
F(h) = 0, 
\]
where $h=1/N$. The function $F$ is defined by 
\[
F(h)= \frac{\Phi(h)}{\sqrt{\Phi(h)^2 + \Psi(h)^2 }} 
-
\frac{\Phi(-h)}{\sqrt{\Phi(-h)^2 + \Psi(-h)^2 }}.
\]
Due to the symmetry $F(-h) = -F(h)$ we have $F(0) = F''(0)=0$. 
Moreover, $F'(0)=0$ because of the fact that $b_1=0$ and so $\Psi(0)=0$. 
The leading order term of Taylor's expansion of $F(h)$ is therefore given as:
\[
F(h)=\frac{1}{6} F'''(0)h^3+O(h^4)
\]
as $h\to 0^+$. Calculation of $F'''(0)$ can be simplified if one uses the substitution 
$\xi(h):=\Psi(h)/\Phi(h)$. 
Notice that $\Phi(h)>0$ for small $|h|\ll 1$ because $\Phi(0)=a_1=g>0$. 
Then $F'''(0)=-6\xi'(0)\xi''(0)$. 
Using (\ref{recurrent}) we finally obtain 
\[
F'''(0)=3\frac{b_2^2}{a_1^3}
\left[ 
\frac{a_2}{2}-\frac{1}{3}\frac{a_1 b_3}{b_2}
\right].
\]
Since $F(h)\equiv 0$ we finally deduce the condition
\[
\frac{a_2}{2} -\frac{1}{3} \frac{a_1 b_3}{b_2} = 0, 
\]
which has to be satisfied by the minimizer $X$ in the limit $h\to 0^+$, 
i.e. for $N\to \infty$. 
Now calculating the corresponding terms $a_1, a_2, b_2$ and $b_3$ 
from the recurrent relations (\ref{recurrent}) we obtain 
\[
a_1=g, 
\quad a_2 = \partial_u g, 
\quad b_2= k g^2, 
\quad b_3= \partial_u(k g^2) + g k \partial_u g.
\] 
Hence
\[
\frac{1}{2}\frac{\partial_u g}{g}-\frac{1}{3}\frac{\partial_u(k g^2)
+g k \partial_u g}{k g^2}=0 
\qquad\hbox{on}\ \ \Gamma=\{\vecx(u);\ u\in[0, 1]\},
\]
which is clearly equivalent to the statement:
$|k|g^{3/2} = constant$, i.e. 
\begin{equation}
|k|^{2/3} g = constant
\label{23condition}
\end{equation}
on the curve $\Gamma=\{\vecx(u);\ u\in[0, 1]\}$. 
It corresponds to the curvature adjusted tangential velocity with 
$\varphi(k)=|k|^{2/3}$ and the extended relative local length (see (\ref{rellength}))
\[
\ratio_{\varphi}(u)=\frac{g(u)}{L}\frac{\varphi(k(u))}{\<\varphi(k(\cdot))\>} \equiv 1, \quad
\hbox{for each}\ u\in[0, 1]. 
\]
In other words, such a curvature adjusted tangential velocity yields 
(in the limit of number of grid points tending to infinity) 
the best possible approximation of the evolving family of planar curves as far as 
the minimization of the error between the length of a curve and the length of its polygonal 
approximation is concerned. 

Figure~\ref{fig:distribution}~(e) displays an example of the length discrepancy minimizing redistribution 
for the case when 
\[
\varphi(k)=|k|^{2/3}.
\] 
The values of length defect $\Delta_L:=1-{\mathcal L}(X)/L$ for several $X$'s corresponding 
the cases in Figure~\ref{fig:distribution} and the length optimal case 
are calculated as in Table~\ref{tbl:length-defect}. 

{\bf The area discrepancy.}\ 
Similarly as in the case of the length functional 
we can ask the question what is an asymptotically optimal redistribution of points on a given curve 
$\Gamma$ such that the discrepancy in areas enclosed by the curve $\Gamma$ and its polygonal approximation is minimal. 

The area $A$ enclosed by the curve $\Gamma$ is given by $A=\frac12 \int_\Gamma \det(\vecx,\partial_s\vecx)\,ds$. 
The area enclosed by the closed polygonal curve with vertices $X=\{\vecx_1, \vecx_2, \cdots, \vecx_N\}$ is given by 
\[
{\mathcal A}(X) = \frac12 \sum_{j=0}^{N-1} \det(\vecx_j,\vecx_{j+1}-\vecx_j), 
\]
where $\vecx_0=\vecx_N$. 
Clearly, ${\mathcal A}(X) = \frac12 \sum_{j=0}^{N-1} \det(\vecx_j,\vecx_{j+1})$. 
The first order necessary condition for a set $X=\{\vecx_1, \vecx_2, \cdots, \vecx_N\}$ of vertices 
to be a minimizer of the area discrepancy functional 
\[
\min_{X\subset\Gamma}(A-{\mathcal A}(X))^2,
\]
where the vertices are constrained to belong to the curve $\Gamma$ reads as follows:
\[
{\mathcal A}^\prime_{\vecx_i}(X) \vecy=0\quad 
\hbox{for any} \ \vecy \sim \vect_i,\quad 
i=1, 2, \cdots, N. 
\]
We have
\[
{\mathcal A}^\prime_{\vecx_i}(X) \vecy
= \frac12 \left(\det(\vecx_{i-1},\vecy) + \det(\vecy, \vecx_{i+1})\right)
=-\frac12 \det(\widetilde{\vect}_i,\vecy), \quad
\widetilde{\vect}_i=\vecx_{i+1}-\vecx_{i-1}
\]
for any direction $\vecy\in {\mathbb R}^2$. 
Hence, the necessary condition for a set $X=\{\vecx_1, \cdots, \vecx_N\}\in {\mathbb R}^{2\times N}$ of points 
belonging to $\Gamma$ to be a minimizer of the functional $X\mapsto (A-{\mathcal A}(X))^2$ can be rewritten as:
\begin{equation}
\widetilde{\vect}_i\ /\!/\ \vect_i
\ \Leftrightarrow\ 
\widetilde{\vect}_i\perp\vecn_i,\quad 
i=1, 2, \cdots, N. 
\label{necessary-condA}
\end{equation}
A graphical description of the necessary condition is depicted in Figure~\ref{fig:necessary-cond}~(b). 

Now, calculating the difference $\widetilde{\vect}_i$ by means of 
the Taylor series and using the facts: $\det(\vect_i,\vecn_i) = 1$ and $\det(\vect_i,\vect_i) = 0$ 
we obtain that condition (\ref{necessary-condA}) can be restated as:
\[
\frac{1}{3} b_3 h^3 + O(h^5) = 0
\]
as $h\to0^+$. 
Therefore, in the limit $h=1/N \to 0$, we conclude that $b_3=0$. 
Since $b_3=\partial_u(k g^2)+g k \partial_u g = k g^2 \partial_u\left( \ln(|k|g^2 ) + \ln g\right)$, 
we obtain $g^3 |k|= constant$ on $\Gamma$, i.e. $g |k|^{1/3} = constant$. 
It means that the area discrepancy curvature adjusted tangential velocity has the shape function 
\[
\varphi(k)=|k|^{1/3}. 
\]

Figure~\ref{fig:distribution}~(f) displays an example of the area discrepancy minimizing redistribution 
for the case when $\varphi(k)=|k|^{1/3}$. 
The values of area defect $\Delta_A:=1-{\mathcal A}(X)/A$ for several $X$'s corresponding 
the cases in Figure~\ref{fig:distribution} and the area optimal case 
are calculated as in Table~\ref{tbl:length-defect}. 
\section{\normalsize Numerical approximation scheme}\label{sec:scheme}

The purpose of this section is to construct a numerical approximation scheme for 
solving the governing equation (\ref{eq:equation_position}) for the position vector and 
the tangential velocity equation (\ref{eq:tangential_velocity}) satisfying 
renormalization constraint (\ref{eq:tangential_velocity_av=0}). 

{\bf Scheme.}\ 
For a given initial $N$-sided polygonal curve $\calP^0=\bigcup_{i=1}^N\calS_i^0$, find a family of $N$-sided polygonal curves $\{\calP^j\}_{j=1, 2, \cdots}$, $\calP^j=\bigcup_{i=1}^N\calS_i^j$, where $\calS_i^j=[\vecx_{i-1}^j, \vecx_{i}^j]$ is the $i$-th edge with 
$\vecx_{0}^j=\vecx_N^j$ for $j=0, 1, 2, \cdots$. 
The initial polygon $\calP^0$ is an approximation of $\Gamma^0$ satisfying 
$\{\vecx_i^0\}_{i=1}^N\subset\calP^0\cap\Gamma^0$, 
and $\calP^j$ is an approximation of $\Gamma^t$ at the time $t=t_j$, 
where $t_j=j\tau$ is the $j$-th discrete time ($j=0, 1, 2, \cdots$) if we use a fixed time increment $\tau>0$, or 
$t_j=\sum_{l=0}^{j-1}\tau_l$ ($j=1, 2, \cdots$; $t_0=0$) if we use 
adaptive time increments $\tau_l>0, l=0, \cdots, j-1 $. The updated curve $\{\calP^{j+1}\}$ is determined from the data $\{\calP^j\}$ at the previous time step by 
using discretization in space and time as follows. 

{\bf Discretization in space.}\ 
Let $\calP=\bigcup_{i=1}^N\calS_i$ be an $N$-sided polygonal curve, where 
$\calS_i=[\vecx_{i-1}, \vecx_{i}]$ is the $i$-th edge and $\vecx_i$ is the $i$-th vertex ($i=1, 2, \cdots, N$; $\vecx_0=\vecx_N$). 
We denote the length of $\calS_i$ by $r_i=|\vecx_{i}-\vecx_{i-1}|$. 
The $i$-th unit tangent vector $\vect_i$ can be defined as 
$\vect_i=(\vecx_{i}-\vecx_{i-1})/r_i$. 
Then the $i$-th unit tangent angle $\nu_i$ is obtained from 
$\vect_i=(\cos\nu_i, \sin\nu_i)^{\mathrm{T}}$ in the following way: 
firstly, from $\vect_1=(t_{11}, t_{21})$, we obtain $\nu_1=2\pi-\arccos(t_{11})$ if $t_{12}<0$; $\nu_1=\arccos(t_{11})$ if $t_{12}\geq 0$. 
Secondly, for $i=1, 2, \cdots, N$  we successively compute $\nu_{i+1}$ from $\nu_{i}$: 
\[
\nu_{i+1}=\left\{\begin{array}{@{}ll}
\nu_i+\arcsin(D) & \mbox{if $I>0$}, \\
\nu_i+\arccos(I) & \mbox{if $D>0$}, \\
\nu_i-\arccos(I) & \mbox{otherwise}, 
\end{array}\right.
\quad
\hbox{where}\ D=\det(\vect_i, \vect_{i+1}), 
\quad
I=\vect_i.\vect_{i+1}.  
\]
Finally, we obtain 
$\nu_0=\nu_1-(\nu_{N+1}-\nu_{N})$ and 
$\nu_{N+2}=\nu_{N+1}+(\nu_2-\nu_1)$. 

In order to derive a discrete numerical scheme, 
we follow the flowing finite volume method adopted for curve evolutionary problems as it was proposed by 
Mikula et al. in~\cite{MikulaS2004b, MikulaSB2010}. 
Let us introduce the ``dual'' volume 
$\calS_i^*=[\vecx_{i}^*, \vecx_{i}]\cup[\vecx_{i}, \vecx_{i+1}^*]$ of $\calS_i$, where $\vecx_{i}^*=(\vecx_{i-1}+\vecx_{i})/2$ 
($i=1, 2, \cdots, N$; $\vecx_{N+1}^*=\vecx_1^*$). 
We define the $i$-th unit tangent angle of $\calS_i^*$ by $\nu_i^*=(\nu_i+\nu_{i+1})/2$. 
The $i$-th curvature $k_i$ has the constant value on $\calS_i$, which 
is obtained from integration of $k=\partial_s\nu$ over $\calS_{i}$ with respect to $s$: 
\[
\int_{\calS_{i}}k\,ds
=k_i\int_{\calS_{i}}\,ds
=k_ir_i, \quad
\int_{\calS_{i}}k\,ds
=\int_{\calS_{i}}\partial_s\nu\,ds
=[\nu]_{\vecx_{i-1}}^{\vecx_{i}}
=\nu_i^*-\nu_{i-1}^*. 
\]
Hereafter, 
$\int_{\calS_{i}}{\sf F}\,ds$ means 
$\int_{s_{i-1}}^{s_{i}}{\sf F}\,ds$ for arc-length $s_i$ satisfying $\vecx_{i}=\vecx(s_i, \cdot)$. 
Thus we have $k_i=(\partial_s\nu^{*})_i$, where 
$(\partial_s{\sf F})_i=({\sf F}_{i}-{\sf F}_{i-1})/r_{i}$. 
The $i$-th curvature $k_i^*$ at $\vecx_i$ can be defined as 
$k_i^*=(k_{i+1}+k_i)/2$ which has the constant value on $\calS_i^*$. 

Next we discretize equation (\ref{eq:tangential_velocity}) for the tangential velocity $\alpha$: 
\[
\partial_s(\varphi\alpha)
=\frac{\<f\>}{\<\varphi\>}\varphi-f+\left(\frac{L}{g}\<\varphi\>-\varphi\right)\omega, 
\]
where $f=(\partial_s^2\beta+k^2\beta)\varphi'(k)-k\beta\varphi(k)$. 
Integrating the above equation over $\calS_i$ yields
\[\begin{split}
\int_{\calS_i}\partial_s(\varphi\alpha)\,ds
&=\big[\varphi(k)\alpha\big]_{\vecx_{i-1}}^{\vecx_{i}} \\\
&=\frac{\<f\>}{\<\varphi\>}\int_{\calS_i}\varphi(k)\,ds
-\int_{\calS_i}f\,ds
+\left(L\<\varphi\>\int_{\calS_i}\frac{1}{g}\,ds-\int_{\calS_i}\varphi(k)\,ds\right)\omega.
\end{split}\]
Then
\[\begin{split}
\psi_i
&=\varphi(k_i^*)\alpha_i-\varphi(k_{i-1}^*)\alpha_{i-1} \\
&=\frac{\<f\>}{\<\varphi\>}\varphi(k_i)r_i
-f_ir_i+\left(\calL\<\varphi\>\frac{1}{N}-\varphi(k_i)r_i\right)\omega, \\
f_i
&=\big((\partial_s(\partial_{s^*}\beta))_i
+k_i^2\beta_i\big)\varphi'(k_i)
-k_i\beta_i\varphi(k_i), 
\end{split}\]
where $\beta_i=\beta(\vecx_i^*, k_i, \vect_i)$ is constant on $\calS_i$, 
\[
(\partial_s(\partial_{s^*}\beta))_i
=\frac{(\partial_{s^*}\beta)_i-(\partial_{s^*}\beta)_{i-1}}{r_i}
=\frac{1}{r_i}\big[\partial_s\beta\big]_{\vecx_{i-1}}^{\vecx_{i}}, 
\quad
(\partial_{s*}{\sf F})_i=\frac{{\sf F}_{i+1}-{\sf F}_{i}}{r_{i}^*}, 
\]
$r_i^*=(r_i+r_{i+1})/2$ is the length of $\calS_i^*$, and 
\[
\calL=\sum_{i=1}^Nr_i
\]
is the total length of $\calP$. 
The averages $\<f\>$ and $\<\varphi\>$ are approximated as: 
\[
\<f\>
=\frac{1}{\calL}\sum_{i=1}^Nf_ir_i, \quad
\<\varphi\>
=\frac{1}{\calL}\sum_{i=1}^N\varphi(k_i)r_i. 
\]

To determine $\{\alpha_i\}_{i=1}^N$, 
we have to take account of the renormalization constraint $\<\varphi\alpha\>=0$. 
It can be discretized as: 
\[
\<\varphi\alpha\>
=\frac{1}{\calL}\sum_{i=1}^N\varphi(k_i^*)\alpha_ir_i^*=0. 
\]
Notice that $\calL=\sum_{i=1}^Nr_i=\sum_{i=1}^Nr_i^*$. 
We define a partial sum of $\{\psi_i\}$ by
\[
\Psi_i=\sum_{l=2}^i\psi_l\quad
(i=2, 3, \cdots, N), \quad
\Psi_1=0. 
\]
With this notation we obtain
\[
\varphi(k_i^*)\alpha_ir_i^*
=\varphi(k_1^*)\alpha_1r_i^*+\Psi_ir_i^*. 
\]
Summing the above terms yields
\[
\sum_{i=1}^N\varphi(k_i^*)\alpha_ir_i^*
=\varphi(k_1^*)\alpha_1\calL
+\sum_{i=1}^N\Psi_ir_i^*, \quad
\calL=\sum_{i=1}^Nr_i^*. 
\]
Hence we obtain 
\[
\alpha_1
=-\frac{1}{\calL\varphi(k_1^*)}\sum_{i=2}^N\Psi_ir_i^*, \quad
\alpha_i
=\frac{1}{\varphi(k_i^*)}
\left(\varphi(k_1^*)\alpha_1+\Psi_i\right)\quad
(i=2, 3, \cdots, N). 
\]

{\bf Discretization in time.}\ 
The semidiscretized (in space) evolution equation (\ref{eq:equation_position}) has the form:
\[
\partial_t\vecx_i
=w_i^*(\partial_{s^*}\vect)_i
+\alpha_i(\partial_{s^*}\vecx^*)_i
+F_i^*\vecn_i^*, \quad
w_i^*=w(\vecx_i, \nu_i^*, k_i^*) \quad
F_i^*=F(\vecx_i, \nu_i^*)
\]
for $i=1, 2, \cdots, N$, 
which follows from integration of $\partial_t\vecx=w\partial_s\vect+\alpha\partial_s\vecx+F\vecn$ 
(see  (\ref{eq:equation_position})) over the volume $\calS_i^*$.  The right-hand side can be discretized  as follows:
\[
\partial_t\vecx_i
=\frac{w_i^*}{r_i^*}\left(
\frac{\vecx_{i+1}-\vecx_{i}}{r_{i+1}}
-\frac{\vecx_{i}-\vecx_{i-1}}{r_{i}}\right)
+\frac{\alpha_i}{2r_i^*}(\vecx_{i+1}-\vecx_{i-1})
+F_i^*\vecn_i^*. 
\]

In our approach, we use the semi-implicit numerical scheme for discretization in time:
\[\begin{split}
&\frac{\vecx_i^{j+1}-\vecx_i^{j}}{\tau}
=a_{-}\vecx_{i-1}^{j+1}
+a_{0}\vecx_{i}^{j+1}
+a_{+}\vecx_{i+1}^{j+1}
+F_i^{*j}\vecn_i^{*j}, \\
&
a_{-}=\frac{b}{r_i^j}-a, \quad
a_{0}=-(a_-+a_+), \quad
a_{+}=\frac{b}{r_{i+1}^j}+a, \quad
a=\frac{\alpha_i^j}{2r_{i}^{*j}}, \quad
b=\frac{w_i^{*j}}{r_i^{*j}}, 
\end{split}\]
where 
${\sf F}_i^j$ is the $i$-th data of an $N$-sided polygonal curve 
$\calP^j\approx\Gamma(t_j)$ for $i=1, 2, \cdots, N$. 
In other words, the following tridiagonal matrix has to be solved 
in order to obtain solution in the new $j+1$ time level: 
\[
-a_{-}\tau\vecx_{i-1}^{j+1}+(1-a_0\tau)\vecx_{i}^{j+1}-a_{+}\tau\vecx_{i+1}^{j+1}
=\vecx_{i}^{j}+F_i^{*j}\vecn_i^{*j}\tau \quad
(i=1, 2, \cdots, N)
\]
subject to periodic boundary conditions $\vecx_{0}^{j+1}=\vecx_{N}^{j+1}$, 
$\vecx_{N+1}^{j+1}=\vecx_{1}^{j+1}$. 
We note that the tridiagonal matrix is strictly diagonally dominant provided that  $\tau$ is small enough. The time step $\tau$ can be also chosen adaptively using the following expression
\begin{equation}\label{eq:adaptive_time_step}
\tau_j=
\frac{r_{\min}^j}{4(1+\lambda)}
\left(\frac{w^{j*}_{\max}}{r_{\min}^j}+\frac{|\alpha^j|_{\max}}{2}\right)^{-1}, 
\end{equation}
where 
$\lambda>0$, 
${\sf F}_{\rm min}=\min_{1\leq i\leq N}{\sf F}_i$, 
${\sf F}_{\rm max}=\max_{1\leq i\leq N}{\sf F}_i$ and 
$|{\sf F}|_{\rm max}=\max_{1\leq i\leq N}|{\sf F}_i|$. 
With this choice of $\tau_j$ the solution $\{\vecx_i^{j+1}\}_{i=1}^N$ exists and it is unique.
\section{\normalsize Computational results}\label{sec:results}

This section is devoted to presentation of various numerical experiments. 
In all figures, for the prescribed $\widehat{\tau}>0$, we plot every $\mu\widehat{\tau}$ discrete 
time step using discrete points representing the evolving curve. 
In every $3\mu\widehat{\tau}$ time step, we plot a polygonal curve connecting those points, 
where $\mu=[[T/\widehat{\tau}]/100]$ ($[x]$ is the integer part of $x$), 
and $T$ is approximation of the final time which computed as follows. 
Let $\calA^{t}$ and $\calL^t$ be the enclosed area and the total length of $\calP^j$ at time $t=t_j$. 
In Figures~\ref{fig:loss_of_convexity},~\ref{fig:venus} and~\ref{fig:ya}, 
$T$ is the smallest discrete time $t$ for which both conditions  $|\calA^{t}/\calA^{t-\tau}-1|<\delta$ 
and $|\calL^{t}/\calL^{t-\tau}-1|<\delta$ are satisfied. In all other figures, 
$T$ is the smallest discrete time $t$ such that $\calA^{t}/\calA^{0}<\delta$. 
In figure captions, we provide the number of grid points $N$, 
$\eps\in(0, 1)$ for $\varphi(k)=1-\eps+\eps\sqrt{1-\eps+\eps k^2}$, 
$\kappa_1$ and $\kappa_2$ for the relaxation function $\omega(t)$. 
In all examples, 
the initial discretization is given by the uniform $N$-division of the $u$-range $[0, 1]$. 

{\bf Initial test curves.}\ 
As an initial test examples we  use a circle, an ellipse and the following initial curves 
$\vecx(u, 0)=(x_1(u), x_2(u))^{\mathrm{T}}$ depicted in Figure~\ref{fig:initial-curves}  and parameterized by
\[
x_1(u)=\cos z, \ x_2(u)=0.7\sin z+\sin x_1+x_3^2, \ 
u\in[0, 1], 
\]
where $x_3=\sin(3z)\sin z$ and $z=2\pi u$ (left), and  
\[
x_1(u)=1.5\cos z, \ 
x_2(u)=1.5(0.6\sin z+0.5x_3^2+0.4\sin x_4+0.1\sin x_5), \ 
u\in[0, 1], 
\]
where $x_3=\sin(3z)\sin z$, $x_4=2x_1^2$, $x_5=3e^{-x_1}$ and $z=2\pi u$ (right) (cf.~\cite{MikulaS2001}). 

{\bf Classical curvature flows.}\ 
According to the classical convexification theory and the asymptotic behavior derived by 
M.~A.~Grayson, M.~Gage and R.~S.~Hamilton, 
in the case where the geometric equation is given by $\beta=k$, 
any embedded curve becomes convex in finite time~\cite{Grayson1987}, 
and any convex curve shrinks to a single point. 
Its asymptotic shape is a circle~\cite{GageH1986}. 
See Figure~\ref{fig:GraysonGageHamilton}. 

{\bf Affine curvature flows.}\ 
Convexification also holds true in the case where the evolution law is the so-called affine scale space 
normal velocity $\beta=k^{1/3}$ (cf.~\cite{SapiroT1994}). 
In this case, convexification was proved and the asymptotic behavior was derived by G.~Sapiro and A.~Tannenbaum. 
They showed that any embedded curve shrinks to a single point with an ellipse as 
the asymptotic shape~\cite{SapiroT1994}. 
See Figure~\ref{fig:affine_evolution}. 

{\bf Experimental order of convergence (EOC).}\ 
In the case of the normal velocity $\beta=k^{1/3}$, 
any ellipse shrinks to a point homothetically. 
It means that an ellipse is a self-similar solution to (\ref{geomrov}) (see Figure~\ref{fig:ellipse}). 
By using the explicit self-similar solution, 
the so-called experimental order of convergence (EOC) can be obtained in the following way. 

When the position vector is described by 
$\vecx(u, t)=\eta(t)\vecz(u)$ with $\eta(0)=1$, 
and when the initial curve is an ellipse $\vecx(u, 0)=\vecz(u)=(a\cos(2\pi u), b\sin(2\pi u))^{\mathrm{T}}$ 
with $a, b>0$, for the curvature we obtain 
$k(u, t)=ab\,\eta(t)^{-1}\zeta(u)^{-3/2}$, 
$\zeta(u)=a^2\sin^2(2\pi u)+b^2\cos^2(2\pi u)$. 
In the case when $\beta=k^{1/3}$, 
$\partial_t\vecx\cdot\vecn=(ab)^{1/3}\eta(t)^{-1/3}\zeta(u)^{-1/2}$ holds. On the other hand, we have 
$\partial_t\vecx\cdot\vecn=\partial_t\eta(t)\,\vecz\cdot\vecn=-ab\,\partial_t\eta(t)\,\zeta(u)^{-1/2}$. 
Hence we obtain the rate $\eta(t)=(1-\frac{4}{3}(ab)^{-2/3}t)^{3/4}$ and the extinction time 
$T=\frac{3}{4}(ab)^{2/3}$, which is determined by $\eta(T)=0$. 

Numerical errors at $t=t_j$ with the number of points $N$ can be defined as 
\[
\err_p^j(N)=\left\{\begin{array}{@{}ll}\disp
\max_{1\leq i\leq N}\left|\frac{(x_{1i}^j)^2}{(a\,\eta(t_j))^2}+\frac{(x_{2i}^j)^2}{(b\,\eta(t_j))^2}-1\right|, & 
\mbox{if $p=\infty$}, \\\disp
\left(\frac{1}{N}\sum_{1\leq i\leq N}
\left|\frac{(x_{1i}^j)^2}{(a\,\eta(t_j))^2}+\frac{(x_{2i}^j)^2}{(b\,\eta(t_j))^2}-1\right|^p\right)^{1/p}, & 
\mbox{if $1\leq p<\infty$}, 
\end{array}\right.
\]
where $\vecx_i^j=(x_{1i}^j, x_{2i}^j)^{\mathrm{T}}$ is the $i$-th grid point. 
Therefore we can define the $L^q\big((0, t_M), L^p(0, 1)\big)$ error such as
\[
E_{p, q}(N)=\left\{\begin{array}{@{}ll}\disp
\max_{1\leq j\leq M}\err_p^j, & \mbox{if $q=\infty$}, \\\disp
\left(\frac{1}{M}\sum_{1\leq j\leq M}\big(\err_p^j\big)^q\right)^{1/q}, & 
\mbox{if $1\leq q<\infty$}, 
\end{array}\right.
\]
where $t_M<T$. 
If $E_{p, q}(N)\approx N^{-\mu}$ holds, 
then $E_{p, q}(N/2)\approx 2^\mu N^{-\mu}$ also holds. 
From these approximations, we may expect $\mu\approx \log_2\big(E_{p, q}(N/2)/E_{p, q}(N)\big)$. 
We therefore use the right hand side as an indicator of the numerical convergence, i.e., the 
EOC from the $L^q\big((0, t_M), L^p(0, 1)\big)$ error is defined as
\[
\EOC_{p, q}(N)=\log_2\frac{E_{p, q}(N/2)}{E_{p, q}(N)}. 
\]

Tables~\ref{table:eps=0}, \ref{table:eps=0.1}, \ref{table:eps=0.5}, \ref{table:eps=0.9} indicate 
the $L^q\big((0, t_M), L^p(0, 1)\big)$ errors $E_{p, q}(N)$ and EOC $\EOC_{p, q}(N)$ 
for $N=2^4, 2^5, 2^6, 2^7, 2^8$, 
$(p, q)\in\{1, 2, \infty\}$, and $\eps=0, 0.1, 0.5, 0.9$. 
Parameters are $\tau=0.1N^{-2}$ and $\kappa_1=\kappa_2=100$, and 
sample data are used at the time $t_j=1.5j/M$ with $j=0, 1, \cdots, M=200$. 
Ellipses with the axes ratio $3:1$ are the same as in Figure~\ref{fig:ellipse}. 
In this case, the extinction time $T=\frac{3}{4}(ab)^{2/3}=1.560\cdots>t_M=1.5$. 
Since the governing system is of the parabolic type, 
we have used the time step $\tau\approx N^{-2}$. 
From Tables~\ref{table:eps=0}, \ref{table:eps=0.1}, \ref{table:eps=0.5}, \ref{table:eps=0.9}, 
we can observe $\EOC\approx 2$ independently from the values $\eps$'s. It means that 
these four tangential redistribution methods seem to be almost equally effective in view of the EOC. 

{\bf The length and the area discrepancy.}\ 
An ellipse is a shrinking self-similar solution to (\ref{geomrov}) 
in the case of the normal velocity $\beta=w(\nu)k$, 
as well as the affine curvature flows $\beta=k^{1/3}$. 
Here the weight is 
\[
w(\nu)=\frac{a^2b^2}{2T(a^2\sin^2\nu+b^2\cos^2\nu)}, 
\]
when the axes ratio of ellipse is $a:b$ and the extinction time is $T>0$. 
Similarly as in the previous subsection EOC, 
one can obtain the position vector of a shrinking ellipse by 
$\vecx(u, t)=\eta(t)\vecz(u)$ with the scaling function $\eta(t)=\sqrt{1-t/T}$ ($\eta(0)=1$, $\eta(T)=0$) and 
the initial ellipse $\vecx(u, 0)=\vecz(u)=(a\cos(2\pi u), b\sin(2\pi u))^{\mathrm{T}}$. 
By using this exact solution, 
we will calculate the numerical test of the length and the area discrepancy in the following way. 

Numerical discrepancy of the length and area defects at the time $t=t_j$ can be defined as
\[
\Delta_L^j=\left|1-\frac{\calL^{t_j}}{L^{t_j}}\right|, \quad
\Delta_A^j=\left|1-\frac{\calA^{t_j}}{A^{t_j}}\right|, 
\]
respectively. 
Here $L^t=\eta(t)L^0$ and $A^t=\eta(t)^2A^0$ are the length and area of $\Gamma^t$, 
and $\calL^t$ and $\calA^t$ are the length and area of $\calP^j$ at $t=t_j$, respectively. 
Therefore we can define the $L^q(0, t_M)$ numerical discrepancy of the length and area as follows
\[
\Delta_{*, q}=\left\{\begin{array}{@{}ll}\disp
\max_{0\leq j<M}\Delta_*^j, & \mbox{if $q=\infty$}, \\\disp
\left(\frac{1}{M}\sum_{0\leq j<M}\big(\Delta_*^j\big)^q\right)^{1/q}, & 
\mbox{if $1\leq q<\infty$}, 
\end{array}\right.
\]
where $*=L, A$ and $t_M\leq T$. 

Table~\ref{table:length-area-discrepancy} indicates 
$\Delta_{L, q}$ and $\Delta_{A, q}$ for $q=1, 2, \infty$. 
Parameters are: $N=100$, $\tau=0.1N^{-2}$ and $\kappa_1=\kappa_2=100$, 
the extinction time is $T=1$, 
sample data are used at the time $t_j=Tj/M$ with $j=0, 1, \cdots, M=200$, and 
the axes ratio of the initial ellipse $\Gamma^0$ are $a=3:b=1$. 
As for the shape functions we chose: 
$\varphi(k)=1-\eps+\eps\sqrt{1-\eps+\eps k^2}$ with $\eps=0, 0.5, 0.9, 1$, 
$\varphi(k)=|k|^{2/3}$ and $\varphi(k)=|k|^{1/3}$. 
In Table~\ref{table:length-area-discrepancy}, 
we can observe that the shape function $\varphi(k)=|k|$ attains the minimum value in 
each $\Delta_{*, q}$ for $*=L, A$ and $q=1, 2, \infty$. 

{\bf Weighted curvature flows.}\ 
Asymptotic behavior of solutions to the weighted curvature flow $\beta=w(\nu)k$ 
is related to self-similar shrinking solutions, which need not be unique. 
For details we refer to the paper by Yagisita~\cite{Yagisita2006} and references therein. 
In Figure~\ref{fig:wcf} we show evolving family of plane curves with $\beta=w(\nu)k$. 

{\bf Forced curvature flows.}\ 
The theory of Grayson, Gage and Hamilton was generalized to the case of an 
anisotropic curvature driven flow by K.~S.~Chou and X.~P.~Zhu in~\cite{ChouZ1999a}. 
They considered the case where the evolution law is given by 
$\beta=w(\nu)k+F(\nu)$ with $w(\nu)=\sigma(\nu)+\partial_\nu^2\sigma(\nu)$, 
where $\sigma$ is a given anisotropy function, 
$\sigma(\nu+\pi)=\sigma(\nu)$, and $F(\nu+\pi)=-F(\nu)$. 
They proved that any embedded curve becomes convex in finite time~\cite{ChouZ1999b}, 
and any convex curve shrinks to a single point 
with the asymptotic shape being a self-similar solution to $\beta=w(\nu)k$, 
i.e., the Wulff shape of $\sigma$~\cite{ChouZ1999a}. 
In Figure~\ref{fig:ChouZhu} we illustrate the convexification theory and 
the asymptotic behavior of evolving plane curves due to K.~S.~Chou and X.~P.~Zhu. 

{\bf Loss of convexity phenomenon.}\ 
In the case when the normal velocity is given by $\beta=k+F(\vecx, \nu)$, 
an initially convex curve may loose its convexity in finite time for a special choice of the external force $F$. 
In Figure~\ref{fig:loss_of_convexity} we present such examples with $F$ 
presented by Nakamura et al. in~\cite{NakamuraMHS1999}. 
It is worth noting that the usual crystalline curvature flow is unable to capture 
this convexity-breaking phenomena. 

{\bf Image segmentation by using a gradient flow.}\ 
Let $\gamma(\vecx)>0$ be an inhomogeneous energy density along the curve $\Gamma^t$. 
If $\gamma$ is differentiable, then the $L^2$ gradient flow of the following energy
\[
E(\Gamma^t)=\int_{\Gamma^t}\gamma(\vecx)\,ds
\]
is realized by the geometric equation (\ref{geomrov}) of the form 
$\beta=\gamma(\vecx)k-\nabla\gamma(\vecx)\cdot\vecn$. 
Such a gradient flow can be successfully utilized in various image segmentation problems. 
For example, let an image intensity function be denoted by $I:\setR^2\supset\Omega\to[0, 1]$. 
We assume that $I(\vecx)$ is a piecewise constant function on each pixel. 
Here $I=0$ ($I=1$) corresponds to the black (white) color and $I\in(0, 1)$ corresponds to a scale of gray colors. 
For simplicity, we assume that our target figures are drawn in white color with the black background. 
Then edges of the image correspond to regions where the gradient $|\nabla I(\vecx)|$ is sufficiently large. 
Let us introduce an auxiliary density function $\gamma(\vecx)=f(|\nabla I(\vecx)|)$ 
where $f$ is a smooth edge detector function such as $f(s)=1/(1+s^2)$ or $f(s)=e^{-s}$. 
Hence the solution curve $\Gamma^t$ of $\beta=\gamma(\vecx)k-\nabla\gamma(\vecx)\cdot\vecn$ 
has the tendency to minimize the energy $E(\Gamma^t)$. In other words, 
it moves toward the edge in the image on which $|\nabla I(\vecx)|$ is large. 
This is a fundamental idea of image segmentation, and 
it has developed to a sophisticated numerical scheme~\cite{MikulaS2004b, MikulaS2006}. 
An example of such a curve evolution converging to an edge in the given images is depicted 
in Figure~\ref{fig:venus}. 

{\bf Image segmentation for sharp images.}\ 
If contrast of the target image is sufficiently high, 
a simpler scheme described in~\cite{BenesKPSTY2008} can be used. 
We consider the geometric flow $\beta=k+F$ and define the forcing term $F(\vecx)$ as follows:
\[
F(\vecx)=F_{\max}-(F_{\max}-F_{\min})\,I(\vecx)
\quad(\vecx\in\Omega),
\]
where $F_{\max}>0$ corresponds to pure black (background) and $F_{\min}<0$
corresponds to pure white (the object to be segmented). 
Maximal and minimal values determine the final shape because in general $1/F$ is equivalent to 
the minimal radius the curve can attain. 
Choosing an appropriate $F$, 
we can make the final shape be rounded, 
or we can prevent the curve from passing through the outline. 
A segmentation of a given sharp image by means of a plane curve evolution is shown in Figure~\ref{fig:ya}. 
\section{\normalsize Conclusions}\label{sec:conlcusion}

In this paper, we proposed and analyzed a new class of tangential velocities 
by which we can control tangential motion and 
grid point redistribution of plane curves that evolve with the normal velocity depending on a general function of 
the curvature, tangential angle and the position vector. 
The curvature adjusted tangential velocity may not only distribute grid points uniformly along the curve 
but also produce a desirable concentration and/or dispersion depending on the curvature. 
We also demonstrated that curvature adjusted tangential redistribution yields the best possible approximation 
of the evolving family of planar curves as far as the minimization of the error between 
the length and area of a stationary curve and that of its polygonal approximation is concerned. 
Numerical experiments based on semi-implicit numerical flowing finite volume method confirmed capability 
of our new method. 

\bigskip\noindent
{\bf Acknowledgments.}\ 
The authors would like to express our gratitude to the anonymous referees for their comments and suggestions. 


\newpage
\setcounter{section}{3}
\setcounter{figure}{0}
\begin{figure}[ht]
\begin{center}
\begin{tabular}{cc}
\scalebox{0.8}{\includegraphics{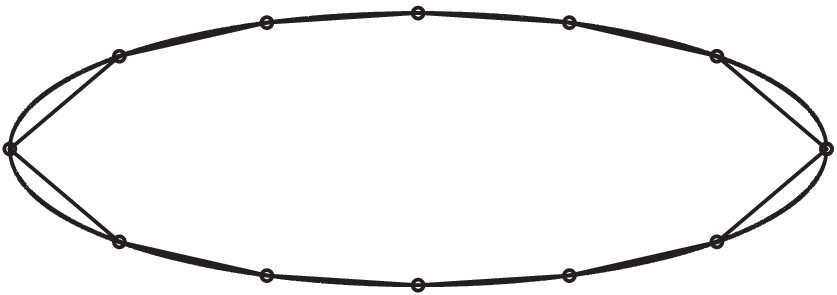}} &\qquad
\scalebox{0.8}{\includegraphics{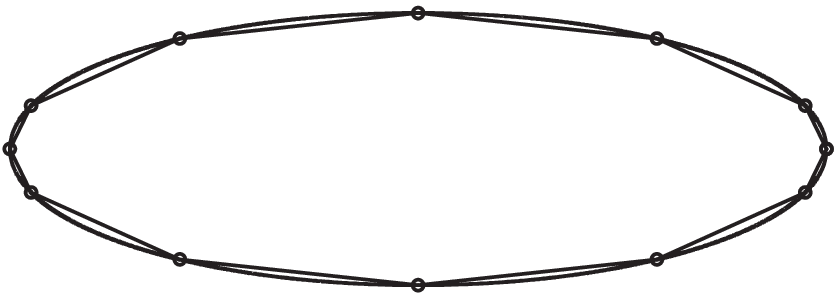}} \\
(a) &\qquad (b) \\[10pt]
\scalebox{0.8}{\includegraphics{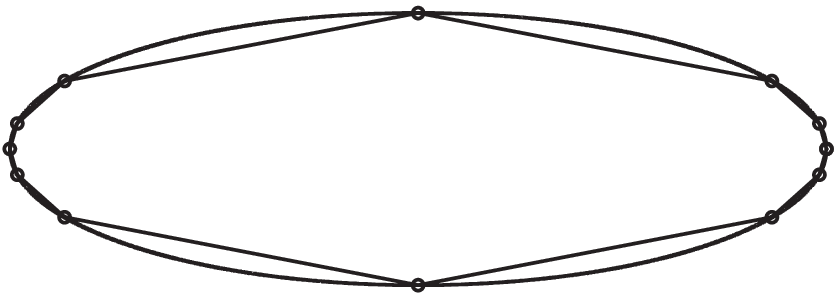}} &\qquad
\scalebox{0.8}{\includegraphics{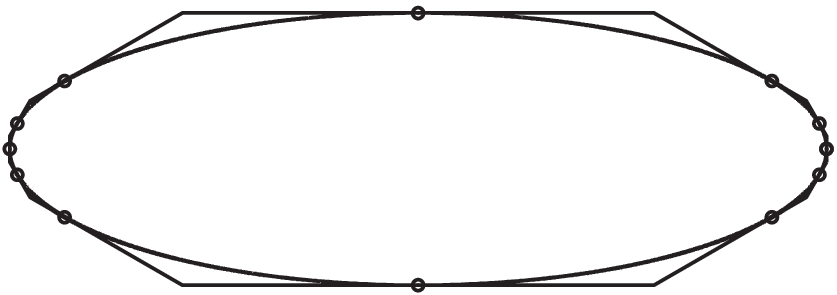}} \\
(c) &\qquad (d) \\[10pt]
\scalebox{0.8}{\includegraphics{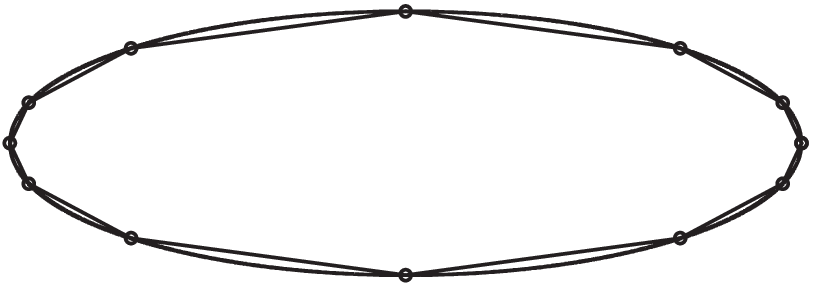}} &\qquad
\scalebox{0.8}{\includegraphics{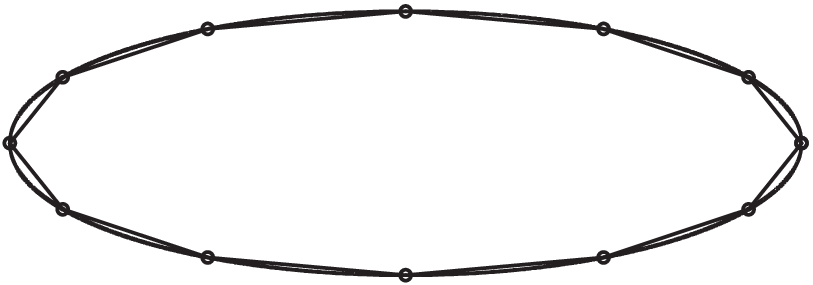}} \\
(e) &\qquad (f)
\end{tabular}
\caption{\small 
Various redistributions of $N=12$ grid points along the ellipse: 
(a)\ uniform redistribution ($\eps=0$), 
(b)\ curvature adjusted redistribution ($\eps=0.9$), 
(c)\ curvature adjusted redistribution ($\eps=1$), 
(d)\ crystalline redistribution (corresponding to $\eps=1$),
(e)\ the length discrepancy minimizing curvature adjusted redistribution with $\varphi(k)=|k|^{2/3}$, and 
(f)\ the area discrepancy minimizing curvature adjusted redistribution with $\varphi(k)=|k|^{1/3}$. 
}
\label{fig:distribution}
\end{center}
\vskip -13pt
\end{figure}
\setcounter{section}{4}
\setcounter{figure}{0}
\setcounter{table}{0}
\begin{figure}[ht]
\begin{center}
\begin{tabular}{cc}
\scalebox{0.25}{\includegraphics{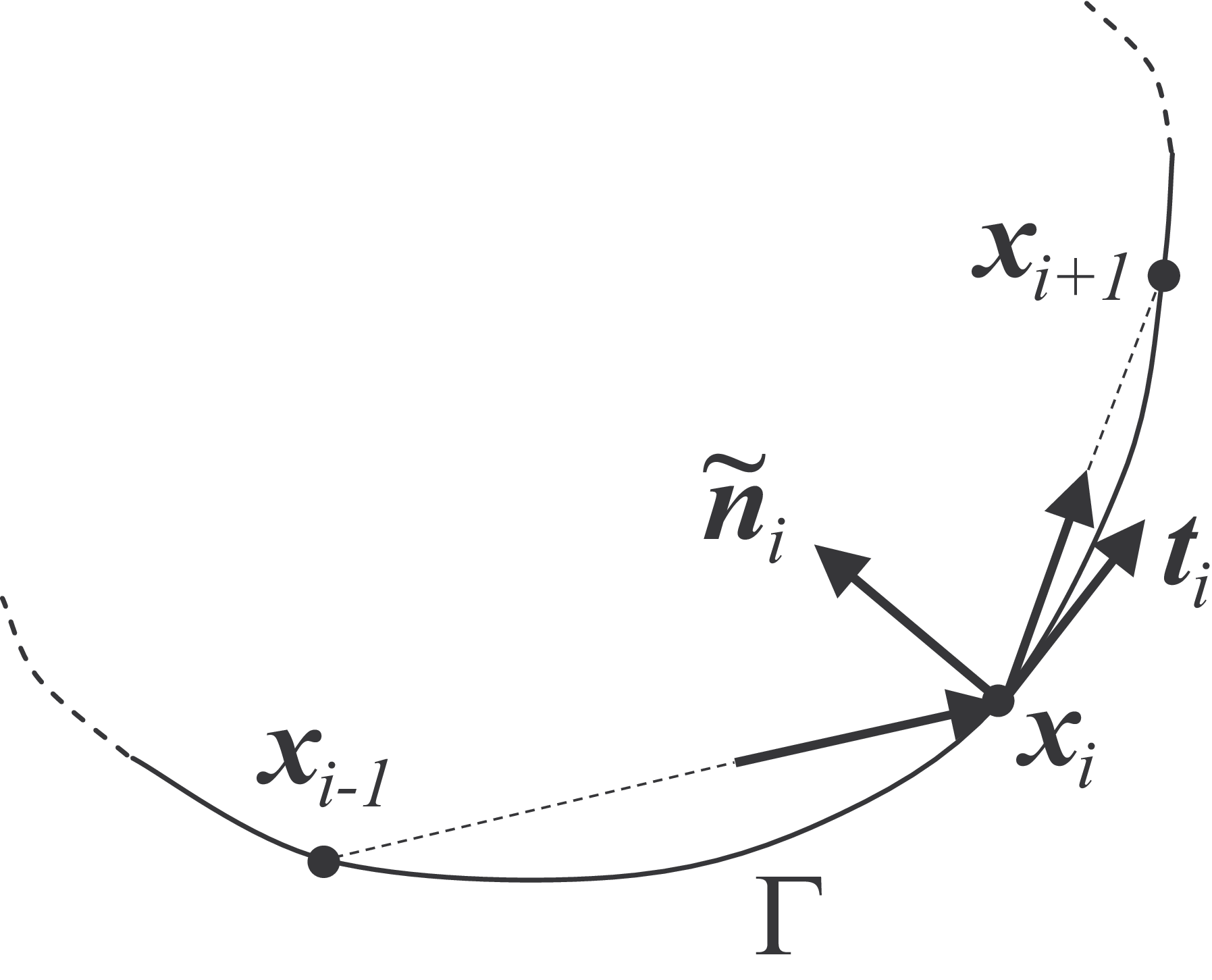}} &\qquad
\scalebox{0.25}{\includegraphics{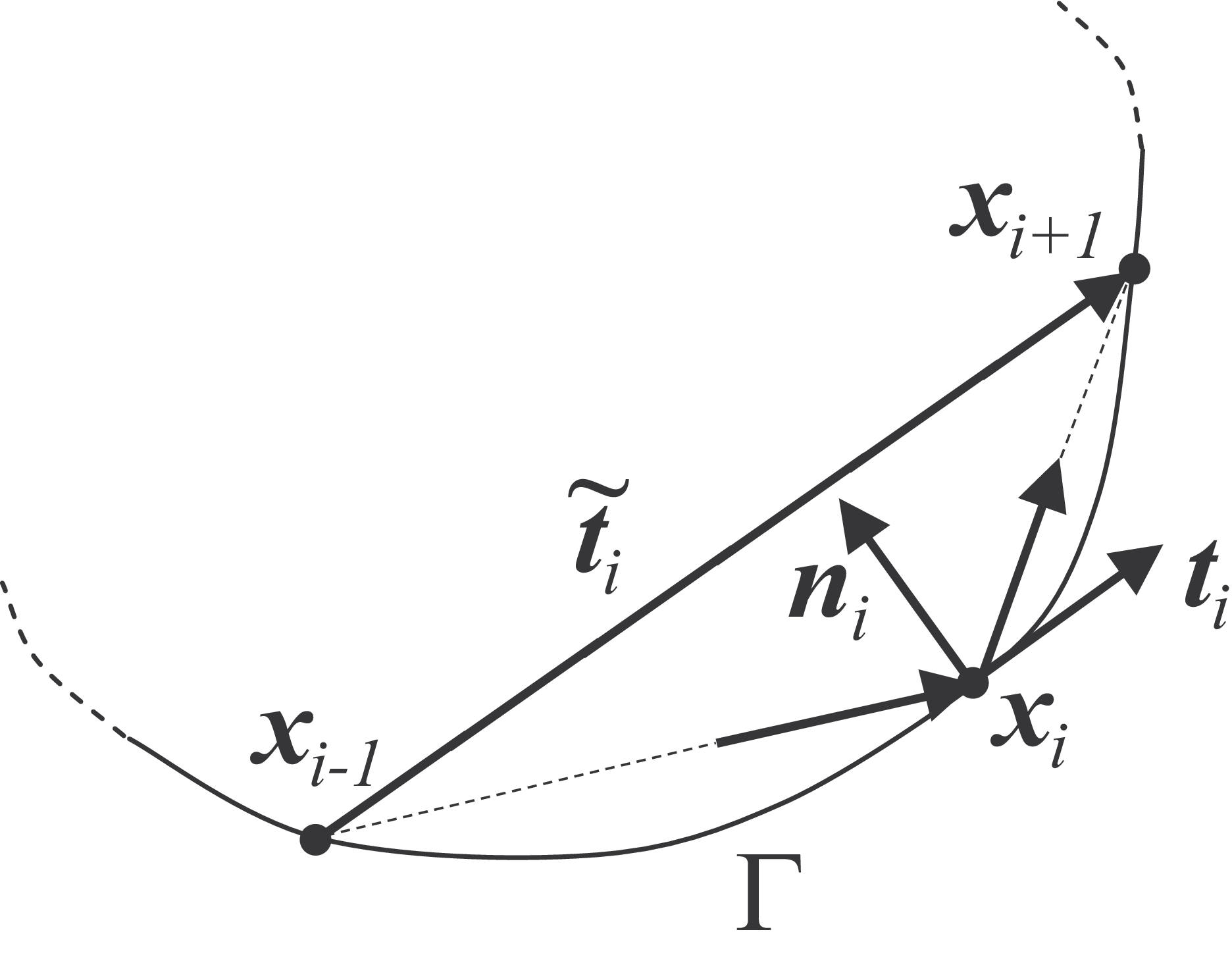}} \\
(a) &\qquad (b)
\end{tabular}
\caption{\small 
A graphical illustrations of the necessary conditions at a point $\vecx_i$ for maximizing of 
the polygonal length discrepancy ${\mathcal L}(X)$ (a) and the polygonal enclosed area ${\mathcal A}(X)$ (b) of a discrete approximation $X=\{\vecx_1, \cdots, \vecx_N\}$ of a smooth curve $\Gamma$, respectively. 
}
\label{fig:necessary-cond}
\end{center}
\vskip -13pt
\end{figure}
\begin{table}[ht]
\begin{center}
\small
\begin{tabular}{|c||c|c|c|c|c|}\hline
& (a) uniform & (b) curv. adj. & (c) crystalline & (e) the length  & (f) the area \\[-10pt]
$X$ & & & & & \\[-10pt]
& $\eps=0$ & $\eps=0.9$ & $\eps=1$ & optimal & optimal \\\hline
$\Delta_L$ & 0.01835 & 0.00789 & 0.00966 & \underline{0.00733} & 0.01085 \\\hline
$\Delta_A$ & 0.05834 & 0.05400 & 0.11998 & 0.06222 & \underline{0.04507} \\\hline
\end{tabular}
\caption{\small 
The values of length defect $\Delta_L$ and area defect $\Delta_A$ 
for the cases in Figure~\ref{fig:distribution}~(a), (b), (c), 
(e) the length optimal case $\varphi(k)=|k|^{2/3}$, and 
(f) the area optimal case $\varphi(k)=|k|^{1/3}$. 
Underlined values are the minimum in each defect. 
}
\label{tbl:length-defect}
\end{center}
\vskip -13pt
\end{table}
\setcounter{section}{6}
\setcounter{figure}{0}
\setcounter{table}{0}
\begin{figure}[ht]
\begin{center}
\begin{tabular}{cc}
\scalebox{0.8}{\includegraphics{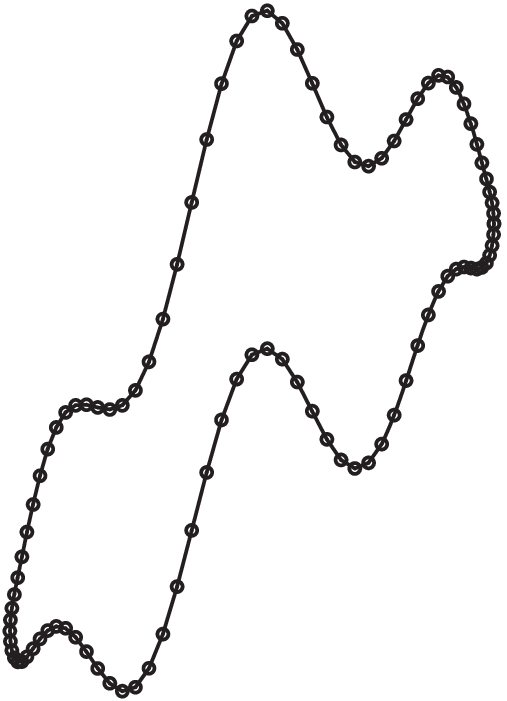}} & \qquad
\scalebox{0.8}{\includegraphics{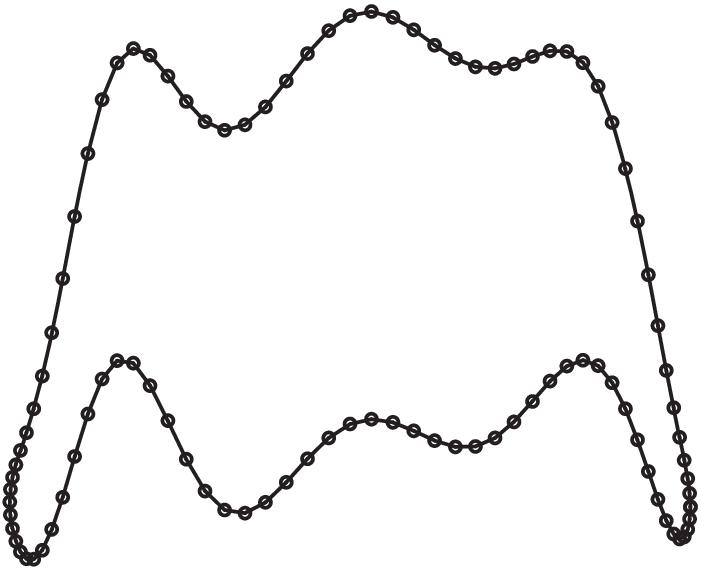}} \qquad
\end{tabular}
\caption{\small 
Initial curves with grid points corresponding to a uniform division of $u\in [0, 1]$ with $N=100$ points.
}
\label{fig:initial-curves}
\end{center}
\vskip -13pt
\end{figure}
\begin{figure}[ht]
\begin{center}
\begin{tabular}{@{}c@{}ccc@{}}
\scalebox{0.8}{\includegraphics{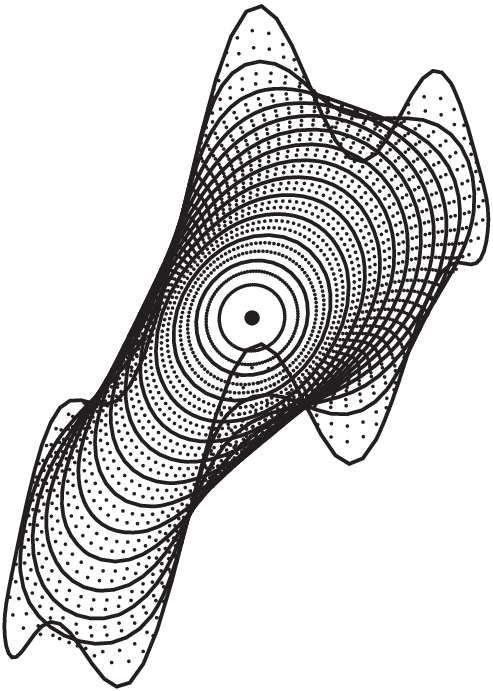}} & 
\scalebox{0.4}{\includegraphics{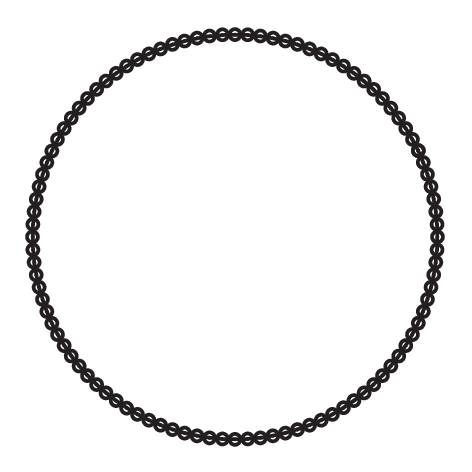}} & 
\scalebox{0.8}{\includegraphics{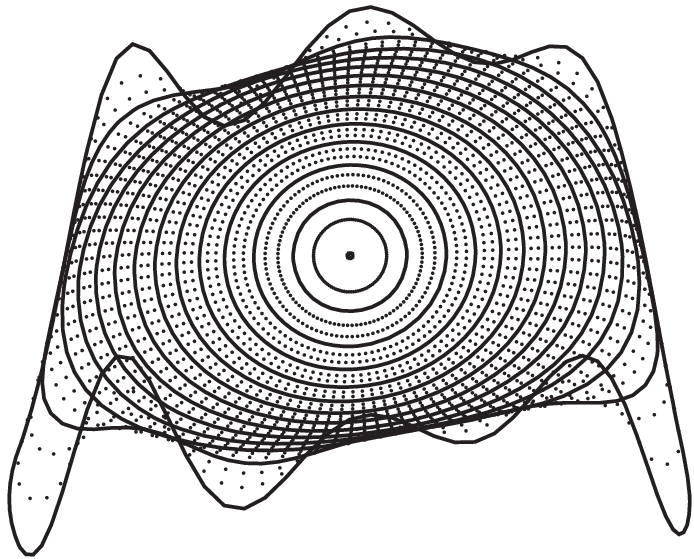}} & 
\scalebox{0.4}{\includegraphics{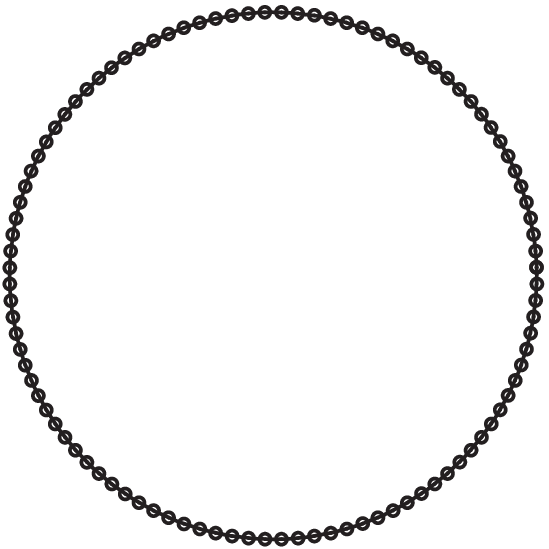}} \\
(a) & (b) & (c) & (d)
\end{tabular}
\caption{\small 
Numerical examples with the normal velocity $\beta=k$. (a) and (c) are evolving curves, respectively. 
(b) and (d) are the corresponding final magnified curves. 
Numerical parameters: $N=100$, $\eps=0.1$, $\kappa_1=\kappa_2=100$, 
$\tau=0.1N^{-2}$,  $\delta=10\tau$,  and $\widehat{\tau}=0.001$. 
}
\label{fig:GraysonGageHamilton}
\end{center}
\vskip -13pt
\end{figure}
\begin{figure}[ht]
\begin{center}
\begin{tabular}{@{}c@{}ccc@{}}
\scalebox{0.8}{\includegraphics{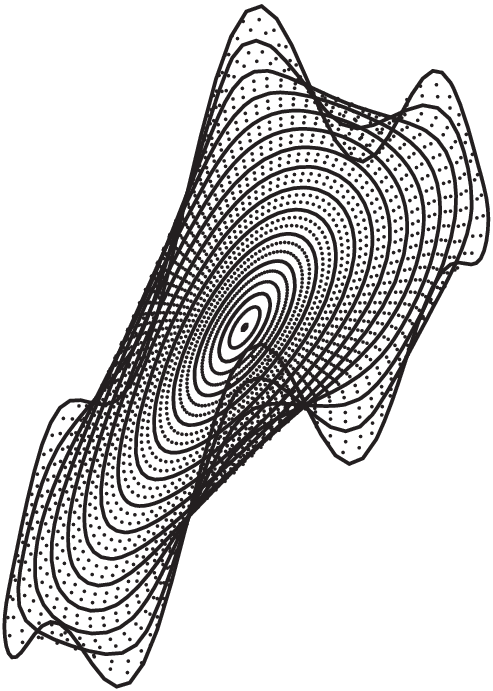}} & 
\scalebox{0.4}{\includegraphics{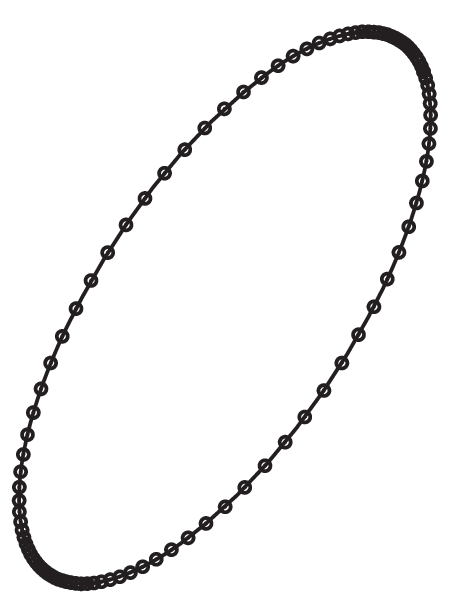}} & 
\scalebox{0.8}{\includegraphics{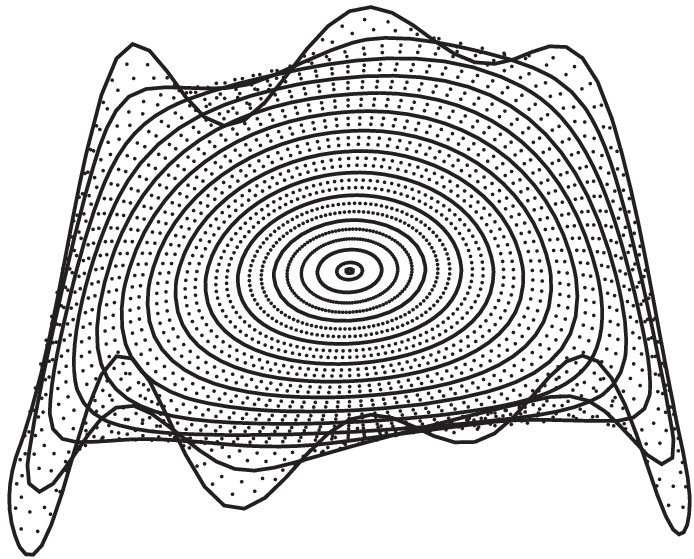}} & 
\scalebox{0.4}{\includegraphics{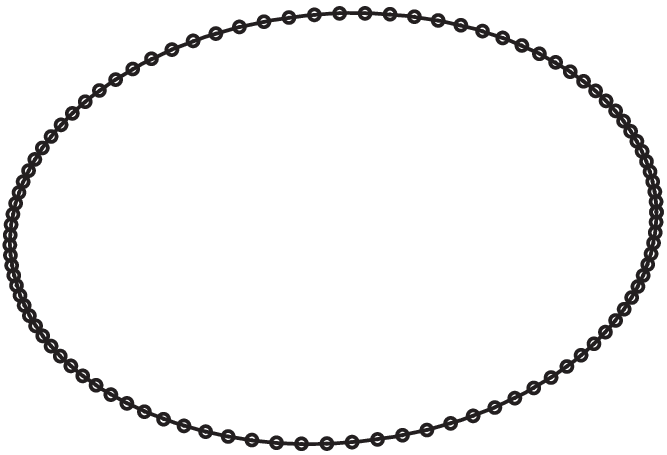}} \\
(a) & (b) & (c) & (d)
\end{tabular}
\caption{\small 
Numerical examples with affine the space scale normal velocity $\beta=k^{1/3}$: 
(a) and (c) are evolving curves, and (b) and (d) are the corresponding final magnified curves. 
In both cases, parameters were chosen: $N=100$, $\eps=0.1$, $\kappa_1=\kappa_2=100$, 
$\tau=0.1N^{-2}$, $\delta=10\tau$, and $\widehat{\tau}=0.001$.
}
\label{fig:affine_evolution}
\end{center}
\vskip -13pt
\end{figure}
\begin{figure}[ht]
\begin{center}
\begin{tabular}{@{}cccc@{}}
\scalebox{0.5}{\includegraphics{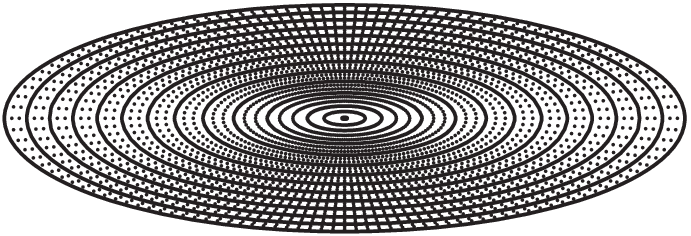}} & 
\scalebox{0.5}{\includegraphics{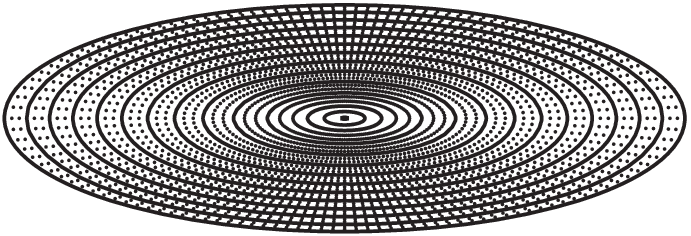}} & 
\scalebox{0.5}{\includegraphics{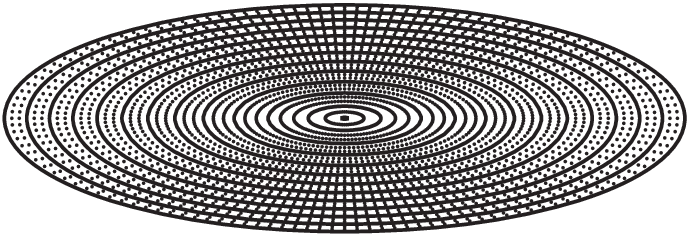}} & 
\scalebox{0.5}{\includegraphics{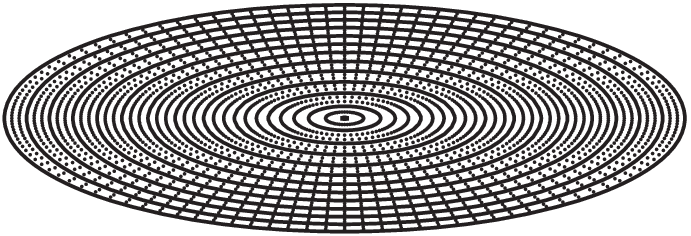}} \\
\scalebox{0.5}{\includegraphics{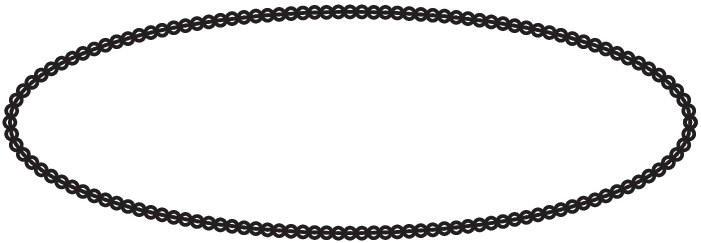}} & 
\scalebox{0.5}{\includegraphics{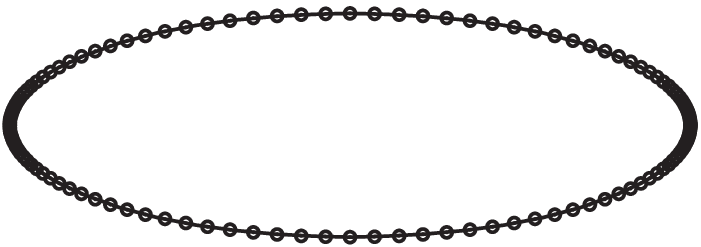}} & 
\scalebox{0.5}{\includegraphics{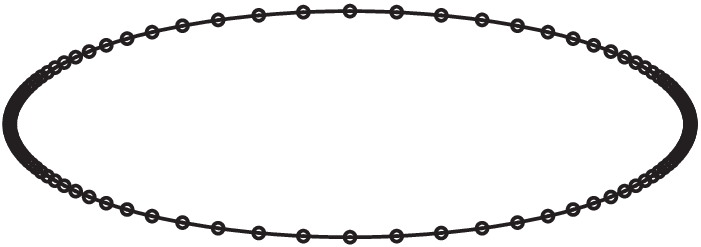}} & 
\scalebox{0.5}{\includegraphics{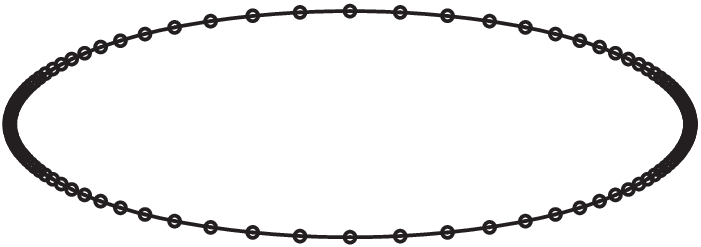}}
\end{tabular}
\caption{\small 
From (far left) to (far right), $\eps=0, 0.1, 0.5, 0.9$. 
Upper figures are evolving ellipses starting from 
the same ellipse with the axes ratio $3:1$. 
Lower figures are the final magnified curves. 
Numerical parameters: $N=128$, $\kappa_1=\kappa_2=100$, 
$\tau=0.1N^{-2}$,  $\delta=10\tau$, 
and $\widehat{\tau}=0.001$.
}
\label{fig:ellipse}
\end{center}
\vskip -13pt
\end{figure}
\begin{table}[ht]
\begin{center}
\small
\begin{tabular}{@{}rccccccc@{}}\hline
$N$ & $p$ & $E_{p, 1}(N)$ & $\EOC_{p, 1}(N)$ & $E_{p, 2}(N)$ & $\EOC_{p, 2}(N)$ & $E_{p, \infty}(N)$ & $\EOC_{p, \infty}(N)$ \\\hline
16 & 1 & 0.5380220 & & 0.5509331 & & 0.6336300 & \\
   & 2 & 1.0534598 & & 1.0704854 & & 1.2142337 & \\
& $\infty$ & 5.9554974 & & 6.0404891 & & 6.8679124 & \\\hline
32 & 1 & 0.1792467 & 1.586 & 0.1825282 & 1.594 & 0.2122894 & 1.578 \\
   & 2 & 0.3427337 & 1.620 & 0.3437584 & 1.639 & 0.3692181 & 1.718 \\
& $\infty$ & 1.8652997 & 1.675 & 1.8668395 & 1.694 & 1.9453397 & 1.820 \\\hline
64 & 1 & 0.0537848 & 1.737 & 0.0552922 & 1.723 & 0.0691816 & 1.618 \\
   & 2 & 0.1018371 & 1.751 & 0.1024454 & 1.747 & 0.1161671 & 1.668 \\
& $\infty$ & 0.5383176 & 1.793 & 0.5391202 & 1.792 & 0.5812710 & 1.743 \\\hline
128 & 1 & 0.0148900 & 1.853 & 0.0154231 & 1.842 & 0.0202166 & 1.775 \\
    & 2 & 0.0281297 & 1.856 & 0.0283963 & 1.851 & 0.0334721 & 1.795 \\
& $\infty$ & 0.1466961 & 1.876 & 0.1470527 & 1.874 & 0.1622977 & 1.841 \\\hline
256 & 1 & 0.0038891 & 1.937 & 0.0040411 & 1.932 & 0.0053923 & 1.907 \\
    & 2 & 0.0073462 & 1.937 & 0.0074271 & 1.935 & 0.0088928 & 1.912 \\
& $\infty$ & 0.0381401 & 1.943 & 0.0382441 & 1.943 & 0.0425560 & 1.931 \\\hline
\end{tabular}

\end{center}
\vskip -13pt
\caption{\small
$\eps=0$}
\label{table:eps=0}
\end{table}
\begin{table}[ht]
\begin{center}
\small
\begin{tabular}{@{}rccccccc@{}}\hline
$N$ & $p$ & $E_{p, 1}(N)$ & $\EOC_{p, 1}(N)$ & $E_{p, 2}(N)$ & $\EOC_{p, 2}(N)$ & $E_{p, \infty}(N)$ & $\EOC_{p, \infty}(N)$ \\\hline
16 & 1 & 0.5290974 & & 0.5419132 & & 0.6255006 & \\
   & 2 & 1.0325342 & & 1.0494491 & & 1.1963774 & \\
& $\infty$ & 5.8214106 & & 5.9054873 & & 6.7533333 & \\\hline
32 & 1 & 0.1743048 & 1.602 & 0.1774370 & 1.611 & 0.2056751 & 1.605 \\
   & 2 & 0.3314933 & 1.639 & 0.3324330 & 1.658 & 0.3553780 & 1.751 \\
& $\infty$ & 1.7951405 & 1.697 & 1.7964082 & 1.717 & 1.8597097 & 1.861 \\\hline
64 & 1 & 0.0518231 & 1.750 & 0.0532112 & 1.738 & 0.0656620 & 1.647 \\
   & 2 & 0.0974615 & 1.766 & 0.0979813 & 1.762 & 0.1092110 & 1.702 \\
& $\infty$ & 0.5117050 & 1.811 & 0.5122756 & 1.810 & 0.5425148 & 1.777 \\\hline
128 & 1 & 0.0142197 & 1.866 & 0.0146984 & 1.856 & 0.0188710 & 1.799 \\
    & 2 & 0.0266349 & 1.872 & 0.0268569 & 1.867 & 0.0309018 & 1.821 \\
& $\infty$ & 0.1377070 & 1.894 & 0.1379516 & 1.893 & 0.1485155 & 1.869 \\\hline
256 & 1 & 0.0036908 & 1.946 & 0.0038251 & 1.942 & 0.0049766 & 1.923 \\
    & 2 & 0.0069038 & 1.948 & 0.0069696 & 1.946 & 0.0081078 & 1.930 \\
& $\infty$ & 0.0354977 & 1.956 & 0.0355656 & 1.956 & 0.0384550 & 1.949 \\\hline
\end{tabular}
\end{center}
\vskip -13pt
\caption{\small
$\eps=0.1$}
\label{table:eps=0.1}
\end{table}
\begin{table}[ht]
\begin{center}
\small
\begin{tabular}{@{}rccccccc@{}}\hline
$N$ & $p$ & $E_{p, 1}(N)$ & $\EOC_{p, 1}(N)$ & $E_{p, 2}(N)$ & $\EOC_{p, 2}(N)$ & $E_{p, \infty}(N)$ & $\EOC_{p, \infty}(N)$ \\\hline
16 & 1 & 0.4042607 & & 0.4129555 & & 0.4910215 & \\
   & 2 & 0.7821867 & & 0.7914808 & & 0.9293693 & \\
& $\infty$ & 4.3834563 & & 4.4198180 & & 5.1595141 & \\\hline
32 & 1 & 0.1310031 & 1.626 & 0.1322110 & 1.643 & 0.1454189 & 1.756 \\
   & 2 & 0.2483055 & 1.655 & 0.2485107 & 1.671 & 0.2563577 & 1.858 \\
& $\infty$ & 1.3333636 & 1.717 & 1.3334130 & 1.729 & 1.3505401 & 1.934 \\\hline
64 & 1 & 0.0362174 & 1.855 & 0.0366554 & 1.851 & 0.0410071 & 1.826 \\
   & 2 & 0.0669561 & 1.891 & 0.0670737 & 1.889 & 0.0698023 & 1.877 \\
& $\infty$ & 0.3448818 & 1.951 & 0.3449364 & 1.951 & 0.3496285 & 1.950 \\\hline
128 & 1 & 0.0093125 & 1.959 & 0.0094400 & 1.957 & 0.0106475 & 1.945 \\
    & 2 & 0.0170431 & 1.974 & 0.0170821 & 1.973 & 0.0178686 & 1.966 \\
& $\infty$ & 0.0863334 & 1.998 & 0.0863494 & 1.998 & 0.0876178 & 1.997 \\\hline
256 & 1 & 0.0023454 & 1.989 & 0.0023787 & 1.989 & 0.0026894 & 1.985 \\
    & 2 & 0.0042800 & 1.994 & 0.0042904 & 1.993 & 0.0044945 & 1.991 \\
& $\infty$ & 0.0215749 & 2.001 & 0.0215790 & 2.001 & 0.0218997 & 2.000 \\\hline
\end{tabular}

\end{center}
\vskip -13pt
\caption{\small
$\eps=0.5$}
\label{table:eps=0.5}
\end{table}
\begin{table}[ht]
\begin{center}
\small
\begin{tabular}{@{}rccccccc@{}}\hline
$N$ & $p$ & $E_{p, 1}(N)$ & $\EOC_{p, 1}(N)$ & $E_{p, 2}(N)$ & $\EOC_{p, 2}(N)$ & $E_{p, \infty}(N)$ & $\EOC_{p, \infty}(N)$ \\\hline
16 & 1 & 0.4408205 & & 0.4417817 & & 0.4758601 & \\
   & 2 & 0.8949908 & & 0.8959769 & & 0.9555951 & \\
& $\infty$ & 5.1442600 & & 5.1492367 & & 5.4758260 & \\\hline
32 & 1 & 0.1161771 & 1.924 & 0.1164446 & 1.924 & 0.1245329 & 1.934 \\
   & 2 & 0.2229378 & 2.005 & 0.2230329 & 2.006 & 0.2323277 & 2.040 \\
& $\infty$ & 1.2060030 & 2.093 & 1.2062437 & 2.094 & 1.2425380 & 2.140 \\\hline
64 & 1 & 0.0290339 & 2.001 & 0.0290913 & 2.001 & 0.0307145 & 2.020 \\
   & 2 & 0.0543111 & 2.037 & 0.0543233 & 2.038 & 0.0558738 & 2.056 \\
& $\infty$ & 0.2827757 & 2.093 & 0.2827892 & 2.093 & 0.2876229 & 2.111 \\\hline
128 & 1 & 0.0072402 & 2.004 & 0.0072539 & 2.004 & 0.0076294 & 2.009 \\
    & 2 & 0.0134384 & 2.015 & 0.0134411 & 2.015 & 0.0137813 & 2.019 \\
& $\infty$ & 0.0690806 & 2.033 & 0.0690829 & 2.033 & 0.0700398 & 2.038 \\\hline
256 & 1 & 0.0018087 & 2.001 & 0.0018121 & 2.001 & 0.0019041 & 2.002 \\
    & 2 & 0.0033502 & 2.004 & 0.0033509 & 2.004 & 0.0034332 & 2.005 \\
& $\infty$ & 0.0171611 & 2.009 & 0.0171616 & 2.009 & 0.0173862 & 2.010 \\\hline
\end{tabular}
\end{center}
\vskip -13pt
\caption{\small
$\eps=0.9$}
\label{table:eps=0.9}
\end{table}
\begin{table}[ht]
\begin{center}
\small
\begin{tabular}{|c||c|c|c||c|c|c|}\hline
$\varphi(k)$ & $\Delta_{L, 1}$ & $\Delta_{L, 2}$ & $\Delta_{L, \infty}$ & $\Delta_{A, 1}$ & $\Delta_{A, 2}$ & $\Delta_{A, \infty}$ \\\hline\hline
$\varphi_0(k)\equiv 1$ &   0.014487 & 0.044087 & 0.461457 & 0.026691 & 0.099960 & 1.119496 \\\hline
$\varphi_{0.1}(k)$ & 0.013191 & 0.039820 & 0.417640 & 0.024203 & 0.089520 & 1.000038 \\\hline
$\varphi_{0.5}(k)$ & 0.005216 & 0.016007 & 0.172409 & 0.009032 & 0.033209 & 0.370852 \\\hline
$\varphi_{0.9}(k)$ & 0.001666 & 0.004468 & 0.047258 & 0.002119 & 0.008126 & 0.092089 \\\hline
$\varphi_1(k)=|k|$ & 0.001051 & 0.001213 & 0.006871 & 0.001291 & 0.002045 & 0.018642 \\\hline
$|k|^{2/3}$ &  0.001124 & 0.002903 & 0.030946 & 0.001427 & 0.005267 & 0.060430 \\\hline
$|k|^{1/3}$ &  0.003680 & 0.011821 & 0.128524 & 0.006159 & 0.023999 & 0.270145 \\\hline
\end{tabular}
\end{center}
\vskip -13pt
\caption{\small 
Numerical discrepancy of the length $\Delta_{L, q}$ and the area $\Delta_{A, q}$ for $q=1, 2, \infty$ 
in various choice of $\varphi(k)$, 
where $\varphi_{\eps}(k)=1-\eps+\eps\sqrt{1-\eps+\eps k^2}$. 
}
\label{table:length-area-discrepancy}
\end{table}
\begin{figure}[ht]
\begin{center}
\begin{tabular}{@{}c@{}ccc@{}}
\scalebox{0.8}{\includegraphics{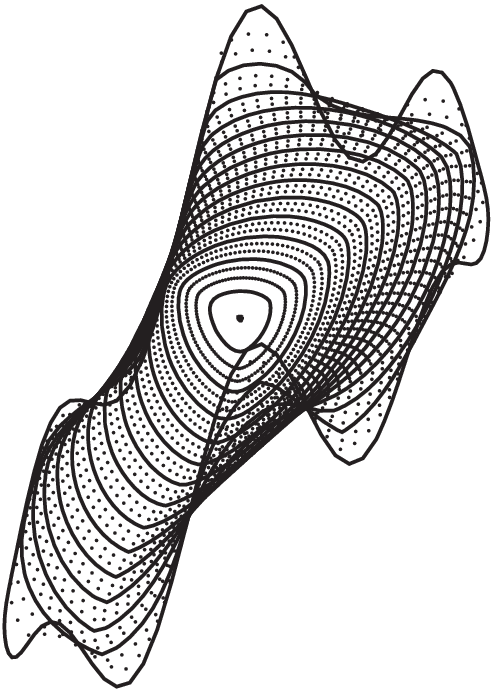}} & 
\scalebox{0.4}{\includegraphics{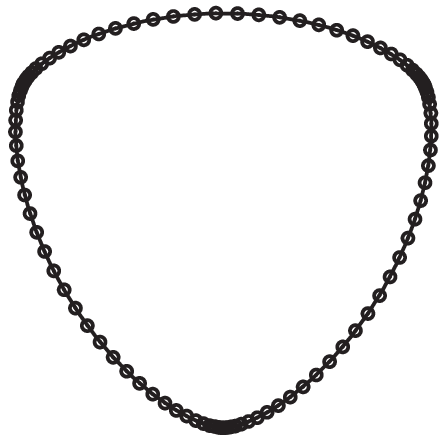}} & 
\scalebox{0.8}{\includegraphics{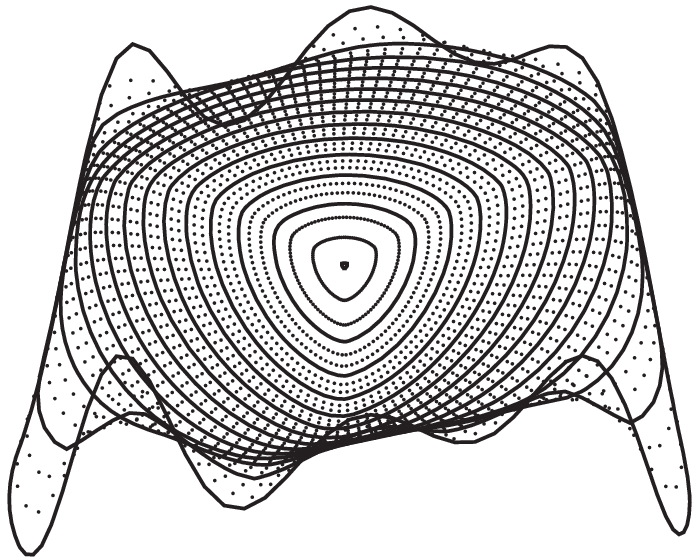}} & 
\scalebox{0.4}{\includegraphics{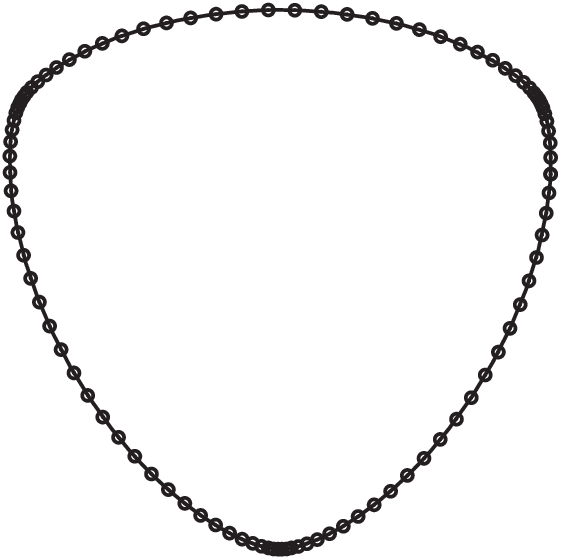}} \\
(a) & (b) & (c) & (d) \\[10pt]
\scalebox{0.8}{\includegraphics{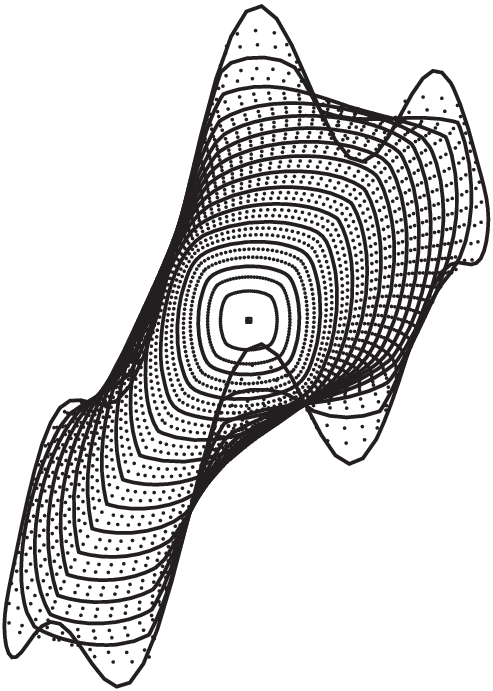}} & 
\scalebox{0.4}{\includegraphics{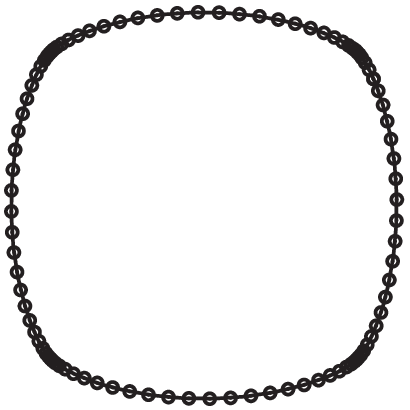}} & 
\scalebox{0.8}{\includegraphics{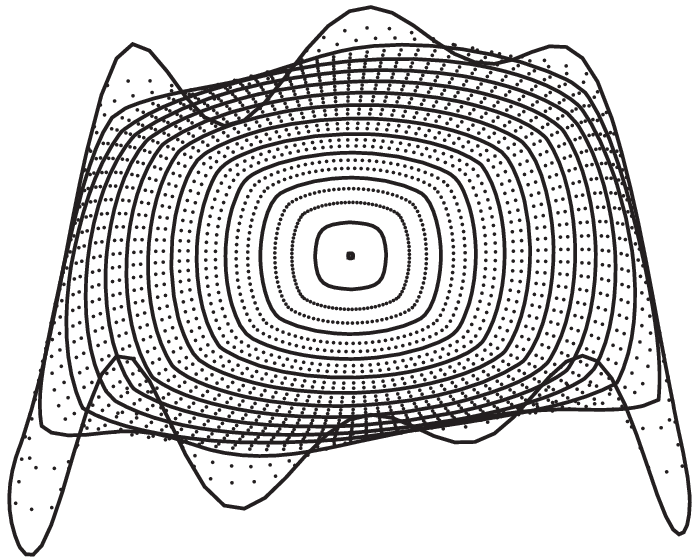}}& 
\scalebox{0.4}{\includegraphics{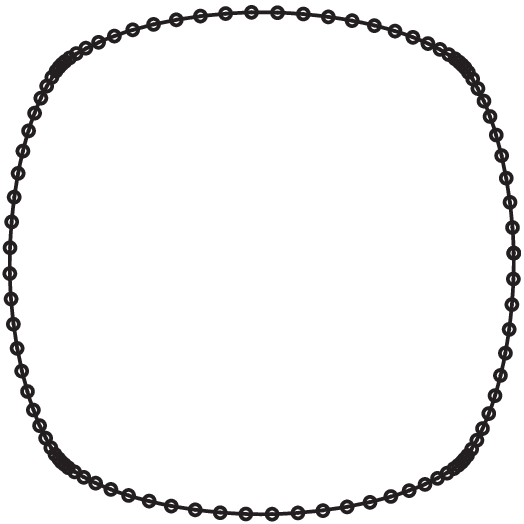}} \\
(e) & (f) & (g) & (h) 
\end{tabular}
\caption{\small 
(a), (b), (c), (d) are in the case where the weight is $w(\nu)=1-(7/9)\cos(3\nu)$, and 
(e), (f), (g), (h) are those where $w(\nu)=1-0.8\cos(4(\nu-\pi/4))$. 
Both weights are used in~\cite{MikulaS2001}. 
(a) and (c) (resp. (e) and (g)) are evolving curves, respectively. 
(b) and (d) (resp. (f) and (h)) are the corresponding final magnified curves. 
In all cases, 
parameters are $N=100$, $\eps=0.1$, $\kappa_1=\kappa_2=100$, 
$\tau=0.1N^{-2}$, 
$\delta=10\tau$, 
and $\widehat{\tau}=0.001$.
}
\label{fig:wcf}
\end{center}
\vskip -13pt
\end{figure}
\begin{figure}[ht]
\begin{center}
\begin{tabular}{@{}c@{}cccc@{}}
\scalebox{0.8}{\includegraphics{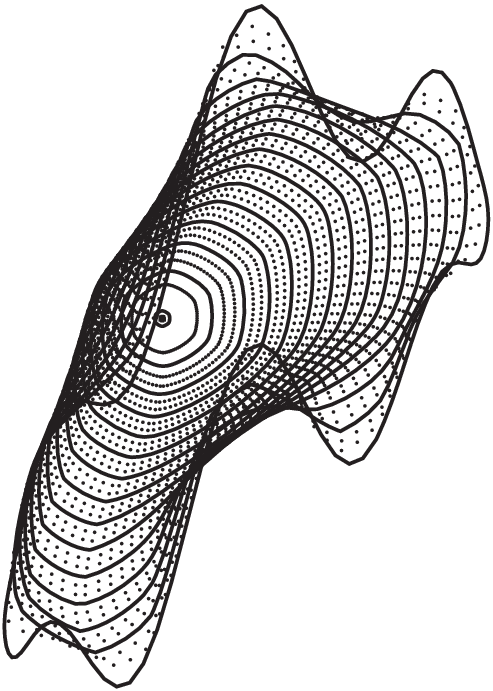}} & 
\scalebox{0.4}{\includegraphics{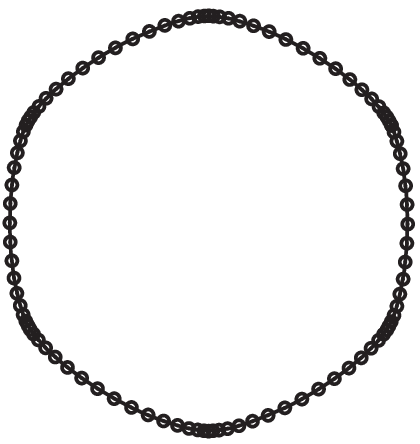}} & 
\scalebox{0.8}{\includegraphics{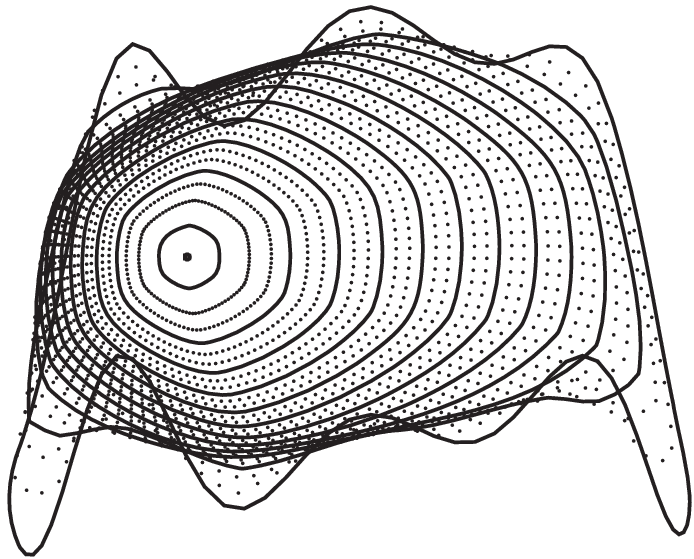}} & 
\scalebox{0.4}{\includegraphics{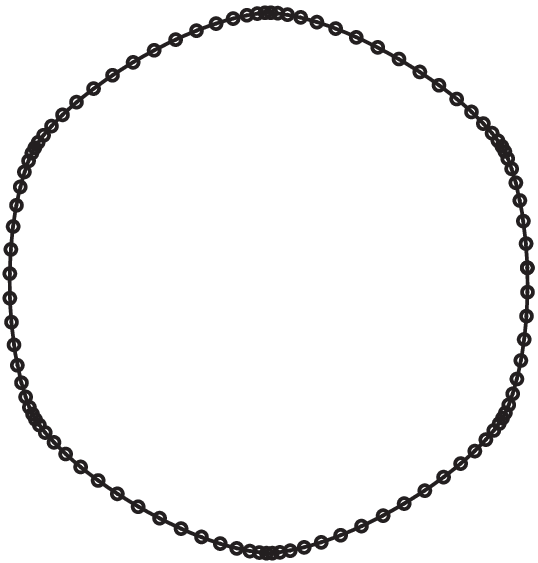}} & 
\scalebox{0.15}{\includegraphics{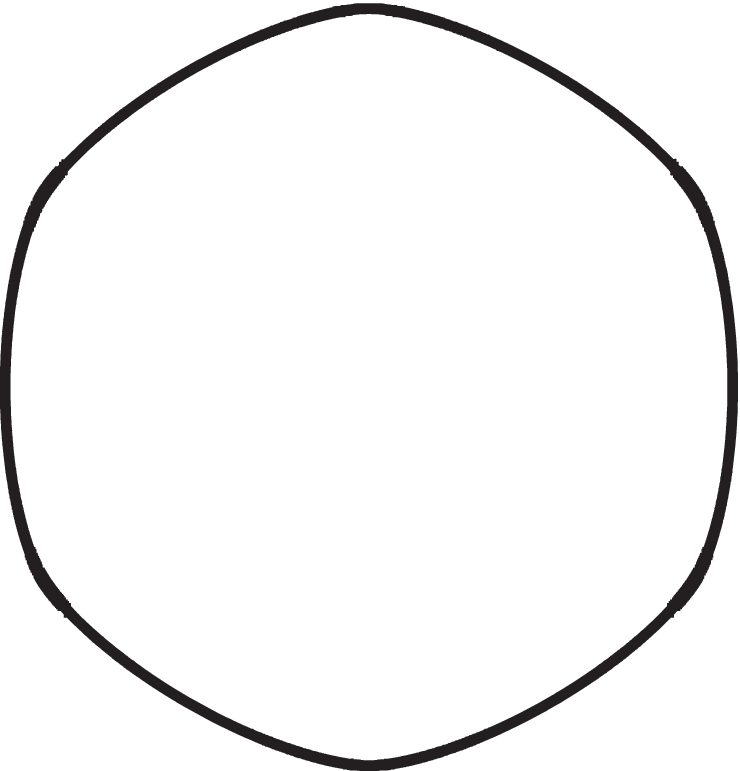}} \\
(a) & (b) & (c) & (d) & (e)
\end{tabular}
\caption{\small 
Numerical examples for the normal velocity $\beta=w(\nu)k+F(\nu)$, 
$w(\nu)=1-0.7\cos(6\nu)$,  $F(\nu)=\sin(\nu)$. 
(a) and (c) are evolving curves, respectively. 
(b) and (d) are the corresponding final magnified curves. 
In both cases, parameters were chosen as:  $N=100$, $\eps=0.1$, $\kappa_1=\kappa_2=100$, 
$\tau=0.1N^{-2}$, 
$\delta=10\tau$, 
and $\widehat{\tau}=0.001$.
The function $\sigma(\nu)=1+(0.7/35)\cos(6\nu)$ is the unique solution of $w=\sigma+\partial_\nu^2\sigma$, 
and (e) is the locus of the boundary of the Wulff shape of $\sigma$, which is given by 
$\sigma(-\vecn)+(\partial_\nu\sigma)\vect$. 
}
\label{fig:ChouZhu}
\end{center}
\vskip -13pt
\end{figure}
\begin{figure}[ht]
\begin{center}
\begin{tabular}{@{}cccc@{}}
\scalebox{0.65}{\includegraphics{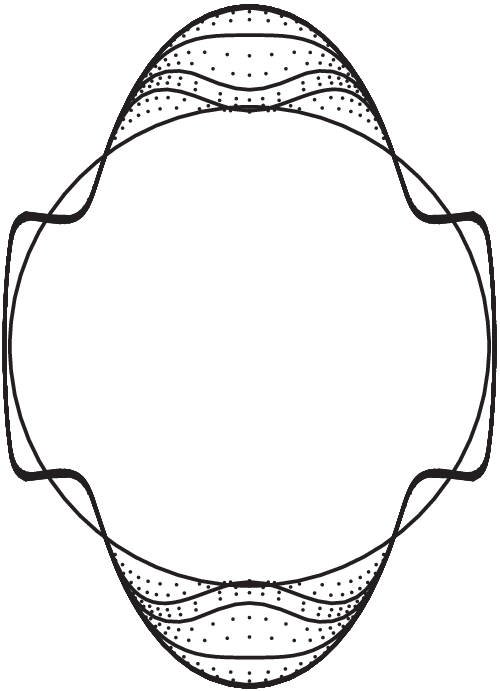}} & \qquad
\scalebox{0.65}{\includegraphics{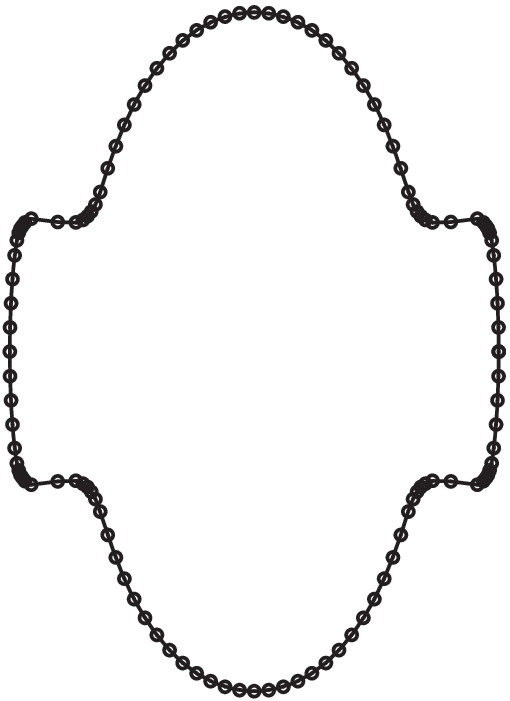}} & \qquad
\scalebox{0.65}{\includegraphics{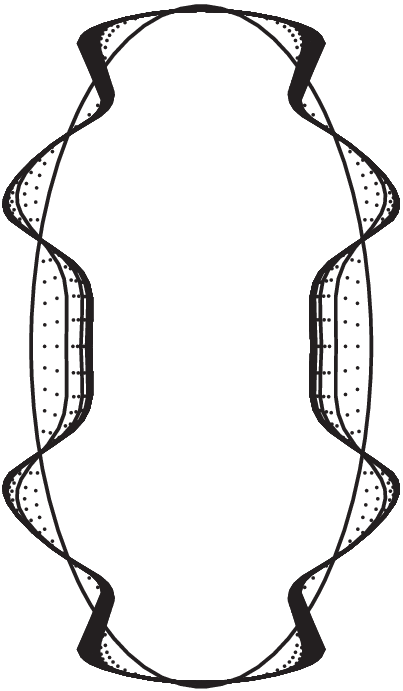}} & \qquad
\scalebox{0.65}{\includegraphics{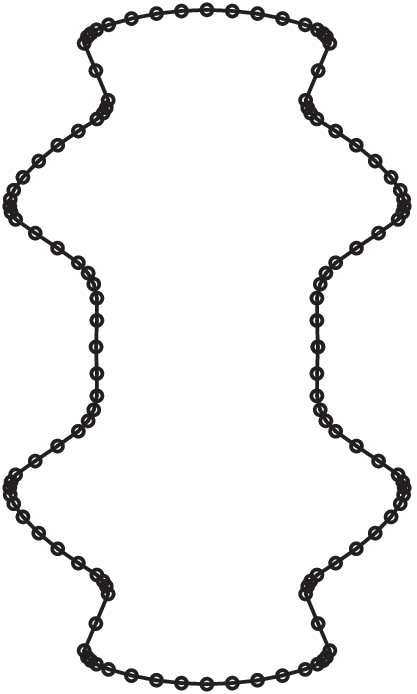}} \qquad \\
(a) & \qquad (b) & \qquad (c) & \qquad (d)
\end{tabular}
\caption{\small 
(a) and (b) correspond to the external force 
$F=-2pq\sin(q(4x_1^2+x_2^2))(-4x_1\sin\nu+x_2\cos\nu)$ with $p=1.25$, $q=3.0$. 
(c) and (d) correspond to the force $F=2pq\pi\cos(q\pi|\vecx|^2)\vecx\cdot\vecn$, 
$p=1.956$, $q=1.15$. 
(a) (resp. (c)) indicates evolving curves starting from 
the unit circle (resp. ellipse with axes ratio $1:2$), and 
(b) (resp. (d)) is the final curve at $T$. 
Parameters are $N=120$, $\eps=0.5$, $\kappa_1=100$, $\kappa_2=0$, 
$\tau=\widehat{\tau}=0.001$, 
and $\delta=10^{-5}$. 
}
\label{fig:loss_of_convexity}
\end{center}
\vskip -13pt
\end{figure}
\begin{figure}[ht]
\begin{center}
\begin{tabular}{@{}ccc@{}}
\scalebox{0.8}{\includegraphics{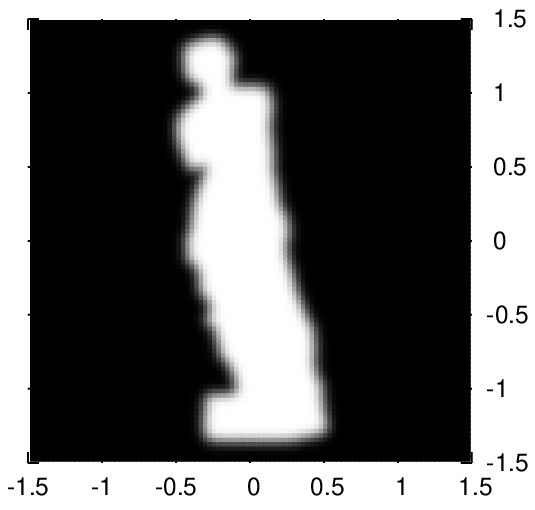}} & 
\scalebox{0.87}{\includegraphics{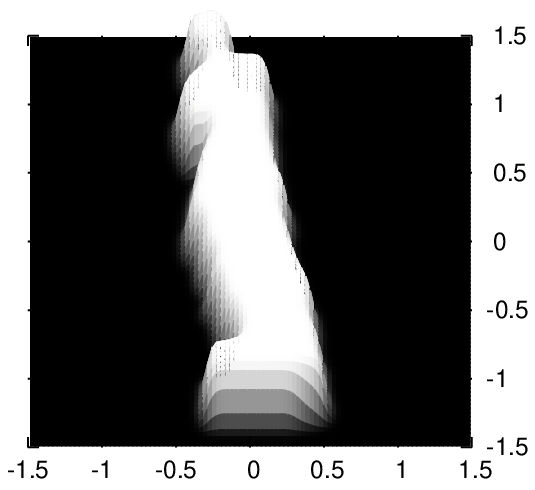}} & 
\scalebox{0.8}{\includegraphics{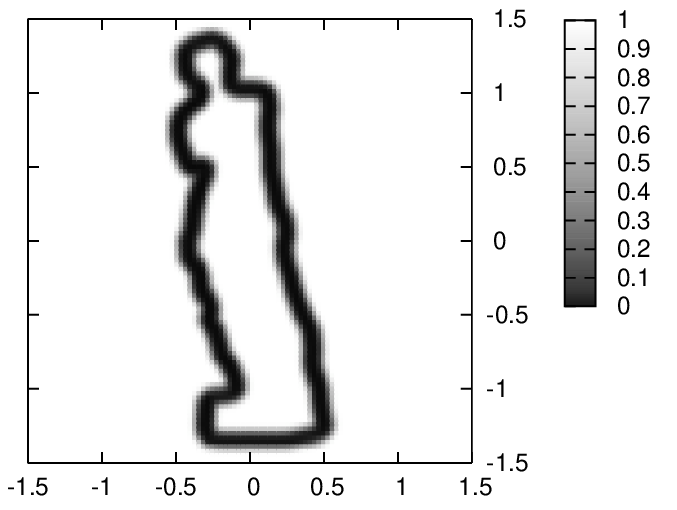}} \\
(a) & (b) & (c)
\end{tabular}\\[10pt]
\begin{tabular}{@{}cc@{}}
\scalebox{0.8}{\includegraphics{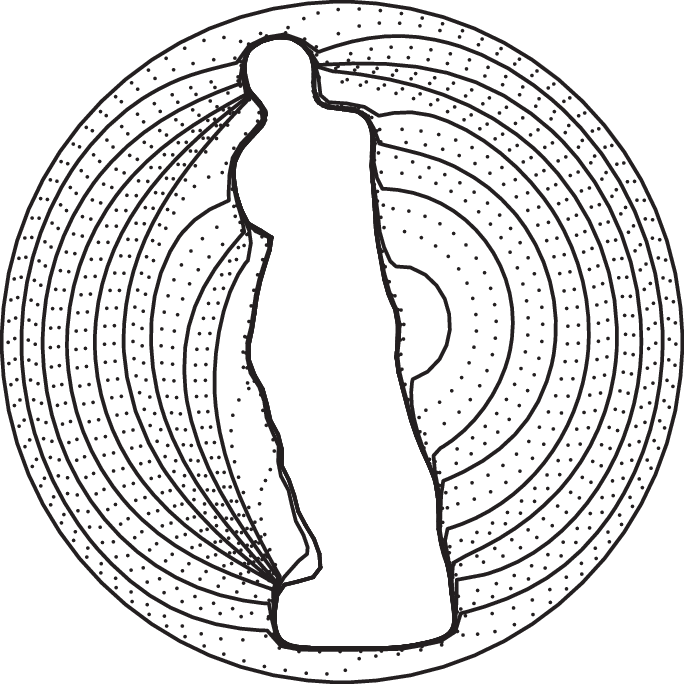}} & \qquad
\scalebox{0.8}{\includegraphics{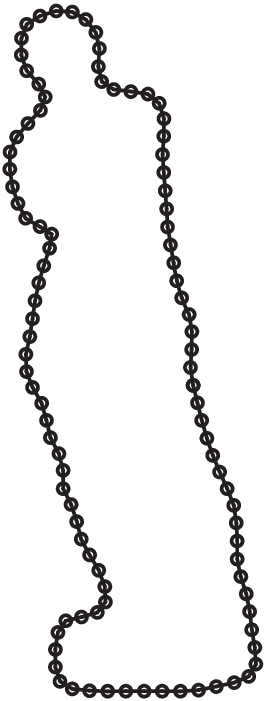}} \\
(e) & \qquad (f)
\end{tabular}
\caption{\small 
(a) indicates the original bitmap image faded shadow ``Venus'' 
in $\Omega=[-1.5, 1.5]^2$ with the 600px resolution, 
(b) indicates image intensity function $I(\vecx)$ in 3D, and 
(c) is the gray scaled auxiliary function $\gamma(\vecx)$. 
(d) indicates evolving curves starting from circle with radius 1.5, and 
(e) is the final curve at $T$. 
Parameters were chosen as: $N=100$, $\eps=0.1$, $\kappa_1=100$, $\kappa_2=0$, 
$\widehat{\tau}=0.001$, 
and $\delta=10^{-5}$. 
Adaptive time step sizing (\ref{eq:adaptive_time_step}) with $\lambda=1$ has been used. 
}
\label{fig:venus}
\end{center}
\vskip -13pt
\end{figure}
\begin{figure}[ht]
\begin{center}
\begin{tabular}{@{}ccc@{}}
\scalebox{0.23}{\includegraphics{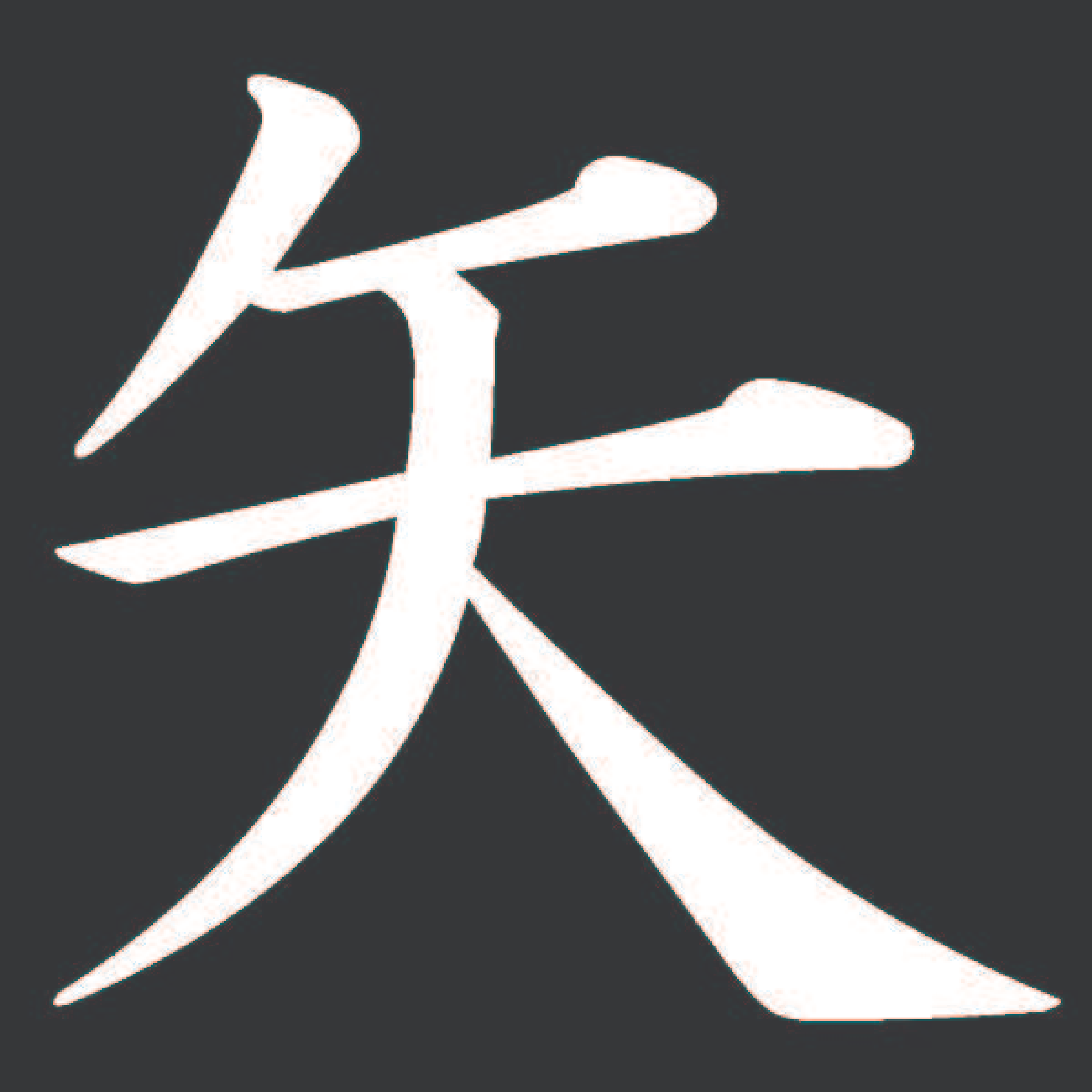}} & 
\scalebox{0.75}{\includegraphics{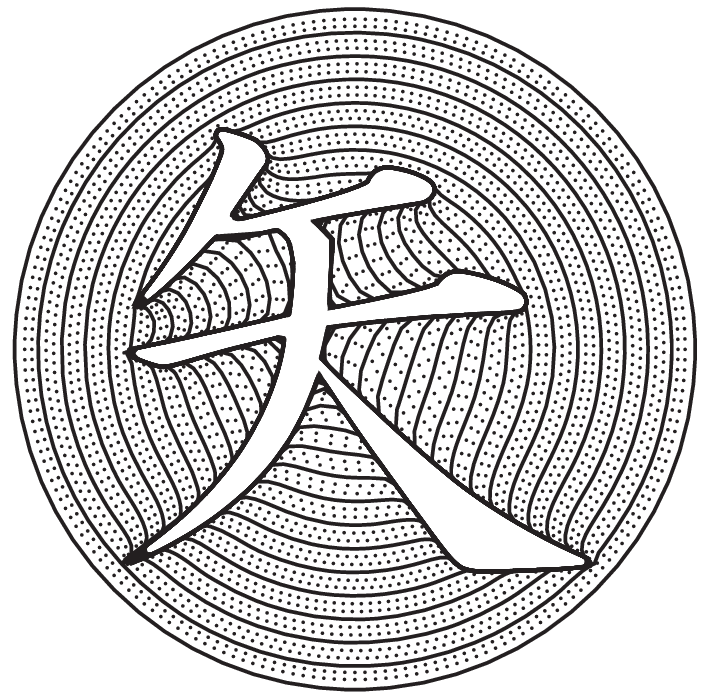}} & 
\scalebox{0.63}{\includegraphics{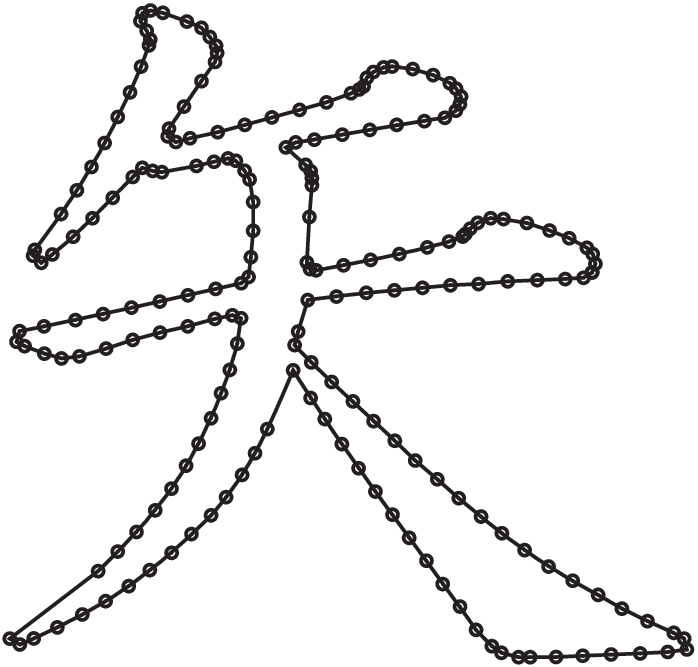}} \\
(a) & (b) & (c)
\end{tabular}
\caption{\small 
(a) indicates the original bitmap image ``ya'' (it is the Chinese character for arrow) in 
$\Omega=[-1.5, 1.5]^2$ with the 600px resolution, 
(b) indicates evolving curves starting from circle with radius 2, and 
(c) is the final curve at $T$. 
Parameters are $N=200$, $\eps=0.1$, $\kappa_1=100$, $\kappa_2=0$, 
$\widehat{\tau}=0.0005$, 
$\delta=10^{-5}/2$, 
$F_{\max}=30$ and $F_{\min}=-30$. 
Adaptive time step sizing (\ref{eq:adaptive_time_step}) with $\lambda=1$ has been used. 
}
\label{fig:ya}
\end{center}
\vskip -13pt
\end{figure}
\end{document}